\documentclass[letterpaper]{amsart}

\usepackage{charter}
\usepackage{amssymb, bm}
\usepackage{amsbsy}
\usepackage{amscd}
\usepackage{amsmath}
\usepackage{amsthm}
\usepackage[
hypertexnames=false,
	hyperindex,
	pagebackref,
	pdftex,
	pdftitle={Algebraic Geometry Fall 2013},
	pdfdisplaydoctitle,
	pdfpagemode=UseNone,
	breaklinks=true,
	extension=pdf,
	bookmarks=false,
	plainpages=false,
	colorlinks,
	linkcolor=blue,
	citecolor=red,
	urlcolor=red,
	pdfmenubar=true,
	pdftoolbar=true,
	pdfpagelabels,
	pdfpagelayout=TwoPage,
	pdfview=Fit,
	pdfstartview=Fit 
]{hyperref}
\usepackage[active]{srcltx}
\allowdisplaybreaks
\numberwithin{equation}{section}

\numberwithin{equation}{section}

\def\BB{{\mathbb B}}

\def\CC{{\mathbb C}}  
\def\DD{{\mathbb D}}

\def\HH{{\mathbb H}}

\def\LL{{\mathbb L}}

\def\PP{{\mathbb P}}
\def\QQ{{\mathbb Q}} 
\def\RR{{\mathbb R}} 
\def\TT{{\mathbb T}} 
\def\XX{{\mathbb X}} 
\def\ZZ{{\mathbb Z}}

\def\pfrak{{\mathfrak{p}}}
\def\ufrak{{\mathfrak{u}}}
\def\vfrak{{\mathfrak{v}}}

\def\gfrak{\mathfrak{g}}

\def\G{\Gamma}
\def\g{\gamma}
\def\alt{{\mathrm alt}} 
\def\bb{{\mathrm bb}}

\def\pr{{\mathrm pr}} 
\def\sst{{\mathrm ss}}
\def\st{{\mathrm st}}

\def\hor{{\mathrm hor}}
\def\ver{{\mathrm vert}}

\def\bs{\backslash}
\def\bbs{\backslash\!\!\backslash}

\def\Gcal{{\mathcal G}} 
\def\Ical{{\mathcal I}} 
\def\Hcal{{\mathcal H}} 
\def\Jcal{{\mathcal J}} 

\def\Lcal{{\mathcal L}}

\def\Ocal{{\mathcal O}}
\def\Pcal{{\mathcal P}}

\def\Tcal{{\mathcal T}}

\def\la{\langle}
\def\ra{\rangle}

\def\pt{{\scriptscriptstyle\bullet}}

\newcommand\Gr{\operatorname{Gr}}
\newcommand\Hom{\operatorname{Hom}}

\newcommand\im{\operatorname{Im}}
\newcommand\iso{\operatorname{Iso}}

\newcommand\spec{\operatorname{Spec}}

\newcommand\proj{\operatorname{Proj}}
\newcommand\sym{\operatorname{Sym}}
\newcommand\GL{\operatorname{GL}}
\newcommand\SL{\operatorname{SL}}
\newcommand\Orth{\operatorname{O}}
\newcommand\SO{\operatorname{SO}}

\newcommand\U{\operatorname{U}}

\newcommand\SU{\operatorname{SU}}

\newcommand\Sp{\operatorname{Sp}}

\newtheorem{theorem}{Theorem}[section]

\theoremstyle{definition}
\newtheorem{definition}[theorem]{Definition}

\newtheorem{example}[theorem]{Example}

\theoremstyle{remark} 
\newtheorem{remark}[theorem]{Remark}

\title[Moduli and locally symmetric varieties]{Moduli spaces  and locally symmetric varieties}
\author[Eduard Looijenga]{Eduard Looijenga}
\address{Mathematical Sciences Center\\  
Jin Chun Yuan West Building, Tsinghua University\\
Haidan District, Beijing 100084 P.R. China}
\email{eduard@math.tsinghua.edu.nl}



\begin{document}

\begin{abstract}
This is a survey paper about moduli spaces that have a natural structure of a (possibly incomplete) locally symmetric variety. We outline the Baily-Borel  compactification for such varieties  and compare it  with the 
compactifications furnished by techniques in algebraic geometry. These differ in general,  but we show that a reconciliation is possible by means of  a generalization of the  Baily-Borel technique for a class of incomplete locally symmetric varieties. 

The emphasis  is here on moduli spaces of varieties other than that of  polarized abelian varieties.
\end{abstract}

\maketitle
\hfill\emph{\begin{small}To Professor Shigeru Mukai on the occasion of his 60th birthday\end{small}}

\section{Introduction}
The goal of this paper is essentially that of the lectures\footnote{Given at the \emph{6th Seasonal Institute of the Mathematical Society of Japan}, Development of Moduli Theory, held in June 2013 in Kyoto}   it is based on, namely to survey and to discuss examples of  locally symmetric (but possibly incomplete) varieties with an interesting modular interpretation, meaning that they can be regarded as moduli spaces of algebraic varieties. In either context, such a  variety  usually comes with a compactification and our focus will eventually  be on their  comparison.

Locally symmetric varieties that are arithmetic quotients of a Siegel half space have such a modular  interpretation naturally, since they can be understood as the moduli spaces of polarized abelian varieties. A fortiori any locally symmetric subvariety of such a variety  has this property as well and it can then almost always be understood as parameterizing polarized abelian varieties with additional structure. We thus obtain an abundance of examples. But not all come that way, at least not in a natural manner, among them the moduli spaces of polarized K3 surfaces. Besides, algebraic geometry often gives us locally symmetric varieties  that are incomplete with respect to their metric. Examples are the moduli spaces of quartic curves, of sextic curves,  of cubic 3-folds and of cubic 4-folds, which, even if we allow them to have innocent singularities miss a locally symmetric hypersurface. In the first two cases the missing locus can be accounted for as parametrizing degenerate objects that become visible only if we  allow the ambient projective plane to degenerate (we thus can account for all genus 3 curves resp.\ for all K3 surfaces of degree 2),  but in the last  two cases it is not clear whether anything like that is true. All these examples have in common that we miss a locally symmetric hypersurface.

This may explain why the compactifications that arise naturally in the two settings, \emph{viz.}\ the ones furnished  by Geometric Invariant Theory and  the Baily-Borel theory, differ. But as we have shown elsewhere \cite{looij1}, \cite{looij2} these two can be reconciled: there is natural extension of the Baily-Borel theory for the kind of locally symmetric varieties that arise here (namely those for which their metric completion  adds a finite union of locally symmetric hypersurfaces) which reproduces the examples that we obtain from Geometric Invariant Theory.

 As in the classical Baily-Borel theory, this comes as a package with topological, analytical and algebraic aspects and has an algebra of (in this case, meromorphic) automorphic forms  as a basic ingredient. A characteristic feature of these compactifications  (which goes back to Satake) is well illustrated by the classical case of the action of a  subgroup $\G\subset \SL (2,\ZZ)$ of finite index on the upper half plane $\HH$: this action is proper and the orbit space $\G\bs \HH$  has in fact the structure of  noncompact Riemann surface  that can be compactified (as a Riemann surface) by filling in a finite set points. The Satake-Baily-Borel approach interchanges the order of things: it tells us first to extend the upper half plane $\HH$ by adding $\QQ\cup\{ \infty\}=\PP^1(\QQ)$  and to equip this union $\widehat\HH$ with a $\SL(2, \QQ)$-invariant topology with the property that (1) each of its points $x$ has a basis of $\G_x$-invariant neighborhoods and (2) the $\G$-orbit space produces a compact surface. This is sensible thing to do from several points of view.  For instance,  $\widehat\HH$ comes with a structure sheaf  $\Ocal_{\widehat\HH}$ of complex valued continuous functions that are analytic on $\HH$. The sheaf of holomorphic differentials on $\HH$ (and its fractional powers) extends to $\widehat\HH$ as an invertible sheaf of $\Ocal_{\widehat\HH}$-modules and what then justifies this set-up is that the $\G$-automorphic forms on $\HH$ can be understood as $\G$-invariant sections of (rational)  powers of this sheaf. 

This is generalized to higher dimensions, but  as one can imagine, things get topologically more involved. The closure of a boundary piece might now have a similar structure so that  the result is a non-locally compact Hausdorff space with a rather intricate topology. This is why  this tends to be a rather technical story. Yet the justification is similar and is perhaps even more compelling in the case of examples borrowed from algebraic geometry, where the reconciliation alluded to above becomes manifest even on the level of algebra: an algebra of invariants with respect to a reductive group  gets identified with  an algebra of meromorphic  automorphic forms. This sometimes allows us to extract structural information on the latter that seems hard to obtain otherwise. We have tried here to slowly lead up to this construction so that the reader will in the end appreciate its naturality.  

In any case, the main characters of our story are the examples, but in order to see them in the light we like to shed on them, we cannot avoid discussing a certain amount  of theory. This is a bit of an equilibrium act, because the underlying theory is notorious for its technicalities and we do not want to get lost in a maze of details. We hope however that once we have presented some examples, the section entitled \emph{Baily-Borel theory and its variations in a nutshell} will serve both the reader who does not want to know these details and the reader who, before delving into the relevant papers, likes to see the general picture first.

\section{Bounded symmetric domains and their classification}

\subsection*{Symmetric spaces}
Let $G$ be a reductive (real) Lie group with compact center. The maximal compact subgroups of $G$ make up a single conjugacy class.  This implies that the set  of maximal compact subgroups of $G$ has the structure of a homogenous manifold (which we denote by $X$): if $K\subset G$ is a maximal compact subgroup of $G$, then $G/K$  has the structure of a  manifold on which $G$ acts by diffeomorphisms  and we have a natural $G$-equivariant identification $G/K\cong X$ with the trivial coset corresponding to the point $o\in X$ associated with $K$. This is the \emph{symmetric space of the noncompact type} attached to $G$. It is known to be diffeomorphic to an open ball.  
Notice that $K$ contains the center of $G$ and so in cases where $X$ is our primary object of study, 
there is no loss in generality in assuming that $G$ has trivial center. This means that $G$ is semisimple and of adjoint type. In particular, $G$ then decomposes as a Lie group naturally  into its simple factors and $X$ decomposes accordingly. Let us therefore assume  in the rest of this section that $G$ is simple, connected and noncompact (we will later allow $G$ to have a finite center).

\subsection*{Invariant complex structures}
To give a complex structure on the tangent bundle of $X$ that is invariant under the $G$-action amounts to giving a complex structure on $T_oX$ (which we suppose handed to us in the form of an action of the circle group  $U(1)$ on $T_oX$) which centralizes the $K$-action. In other words, the complexification $\CC\otimes T_oX$ decomposes 
into a pair of complex conjugate eigenspaces $T^{1,0}_oX$ and $T^{0,1}_oX$ on which  $z\in U(1)$ acts by scalar multiplication by $z$ resp.\ $\bar z$. The condition is then simply that $K$ respects this decomposition.
One  finds that the complex structure on the tangent bundle is integrable so that  $X$ becomes  a complex manifold on which $G$ acts by holomorphic  $X$ automorphisms.  It turns out that this circle group $U(1)$ lies in the image of $K$ and that $K$ is the centralizer of this circle group. Since $-1\in U(1)$ will act as minus the identity on $T_oX$ and centralizes $K$, it is the \emph{Cartan involution of the pair $(G,K)$}: $K$ is also the centralizer of this element. So we here have a distinguished conjugacy class of circle subgroups of $G$ and an invariant complex structure is obtained by  having  these circles parametrized: we then assume given a distinguished conjugacy class $\XX$ of embeddings  $u: U(1)\hookrightarrow G$. We will think of $\XX$ as $X$  endowed with a complex structure. The conjugacy class of $u^{-1}: U(1)\hookrightarrow G$ then gives the  complex conjugate structure on $X$.

\subsection*{An outline of the classification}
This circle group plays also an important role for the classification of the pairs $(G,\XX)$. We outline how this is accomplished. Choose a Cartan subgroup $H$ for $K$. This subgroup is unique up to $K$-conjugacy and will contain $U(1)$. We need two facts, the first one being that $H$ will then also be a  Cartan subgroup for $G$ (so $G$ and $K$ have equal rank). The characters $R\subset \Hom (H, \CC^\times)$ of the adjoint action of  $H$ on $\gfrak_\CC$  make up an irreducible  root system. Since $G$ has trivial center, $R$ spans the abelian group $\Hom (H, \CC^\times)$.  The second fact we need is  that $u: U(1)\hookrightarrow H$ has the remarkable property that its composite with any root is either trivial, the natural inclusion  $U(1)\subset \CC^\times$ or the complex conjugate of the  latter. 

We express this as a property the root system.  We first get rid of $H$ by viewing $R$ as an abstract root system in a real vector space $V$ (in our case that would be $\Hom (H, \CC^\times)_\RR$). For the compact torus  $H$ can then be simply recovered from $R$, or rather from the lattice $Q(R)$ spanned by $R$ in $V$,  as $\Hom (Q(R), U(1))$. We can now represent $u$ as an element $\varpi\in V^*$  via the identity $\alpha (u(z))=z^{\varpi (\alpha)}$, so that  $\varpi |R$ takes values  in $\{ -1,0,1\}$  (it is a so-called a \emph{minuscule coweight}). Then the isomorphism type of $(G, \XX)$ is completely given by the isomorphism type of the pair $(R, \varpi)$. 

This converts the classification into one involving Dynkin diagrams only.  We choose a system of simple roots $(\alpha_1,\dots , \alpha_r)$ of $R$ (where $r$ is the rank of $G$) such that $\varpi$ is dominant: $\varpi (\alpha_i)\ge 0$ for all $i$. If
$\tilde\alpha=n_1\alpha_1+\cdots +n_r\alpha_r$ is the highest root (the coefficients are positive integers), then we must have  $\varpi (\tilde\alpha)=1$, which of course can only happen if there exists a simple root $\alpha_i$ such that $n_i=1$ and $\varpi (\alpha_j)=\delta_{ij}$.  This is what is called a \emph{special vertex} of the Dynkin diagram.  We thus have obtained a pair consisting of a connected Dynkin diagram and a special vertex of it,  given only up to isomorphism in the sense that two special vertices of the same Dynkin diagram $D$ gives rise to isomorphic pairs $(G,\XX)$ if there is an automorphism of $D$ that carries one onto the other. Note that removal of the special vertex yields the Dynkin diagram of $K$.

This also works in the opposite direction in the sense that a pair $(R, \varpi)$ produces $(G, \XX)$ up to isomorphism: the Chevalley construction yields a complex torus $\Hcal:=\Hom (Q(R),\CC^\times)$ and  a complex algebraic group $\Gcal$ containing $\Hcal$ as a Cartan subgroup  and having $R$ as its root system. This construction comes with a natural compact  real form: an anti-involution $c$ of $\Gcal$ which on $\Hcal$ is given by $v\mapsto 1/\bar v$ and has the property that its fixed point subgroup $\Gcal^c$ is compact. Note that $H:=\Hcal^c$ equals $\Hom (Q(R),U(1))$ and that 
$\varpi$ defines  a one-parameter subgroup $u: U(1)\hookrightarrow H$. We then let $G$ be the fixed point subgroup of the anti-involution  of $\Gcal$ defined by $g\mapsto u(-1)c(g)u(-1)^{-1}$. This real form of $\Gcal$ contains $K:=\Gcal^c\cap G$ as a maximal compact subgroup. It is clear
that $K\supset H$. Then $u$ endows $X:=G/K$ with a $G$-invariant structure of  a complex manifold.
\\

With the help of the tables we can now list the cases.  Below we give for the infinite series (the so-called classical cases) the Cartan notation for the associated domains. 

\begin{description}
\item[$A_r$, $r\ge1$] Any vertex of $A_r$ will do, but this is really up to the natural involution of the Dynkin diagram: let us denote by $A_r^{(p)}$, $p=1,2,\dots ,\lfloor  r/2\rfloor$, the case where the special vertex is at position $i$ from an end vertex.  This is Cartan's I$_{p,r+1-p}$.
\item[$B_r$, $r\ge 2$] The end vertex corresponding a long root. We get Cartan's IV$_{2r-1}$.
\item[$C_r$, $r\ge 2$] The unique vertex corresponding a long root (it is an end vertex).  In the Cartan
notation: III$_r$.
\item[$D_r$,  $r\ge 4$]  The special vertices are the three  end vertices. For $r>4$ we have two orbits: 
one consists of the end vertex at a long arm and which  we shall denote $D'_r$ (Cartan IV$_{2r-2}$) and the other  consists of the two end vertices at a short arm, denoted here by $D''_r$ (Cartan II$_{r}$). For $r=4$, the
three end vertices are all equivalent (and indeed,  IV$_6$ and II$_4$ are isomorphic).
\item[$E_6$] The two end vertices at the long arms; they are equivalent under the nontrivial involution of the diagram.
\item[$E_7$] The end vertex on the longest arm.
\end{description}

\section{Classical symmetric domains and Hodge theory}
For $(G, \XX)$ as in the previous section, the complex manifold $\XX$ has the structure of a bounded domain. Harish-Chandra gave a uniform recipe for realizing $\XX$ as a bounded open subset in the  tangent space $T_o\XX$. We will not review this, but we will make his construction explicit for the classical examples. We keep in this section the assumption that $G$ is \emph{simple as a Lie group}, but we no longer wish to assume that $G$ has trivial center. The reason is that we want to consider certain vector bundles over $\XX$ with  $G$-action (lifting of course the given action on $\XX$) and then allow the possibility that the center acts nontrivially on the fibers. Here we are mainly concerned with line bundles that are roots of the  canonical line bundle $\omega_{\XX}$ of $\XX$ ($\omega_{\XX}$ is the sheaf of holomorphic sections of the top exterior power of the holomorphic cotangent bundle of $\XX$ and so comes  naturally with a $G$-action). By  an $N$th `root' of $\omega_{\XX}$ we mean a holomorphic line bundle $\Lcal_\XX$ on $\XX$ with $G$-action such that $\Lcal_\XX^N$ is $G$-equivariantly isomorphic to  $\omega_{\XX}$. The constructions discussed here often lead a natural  (positive) root of  $\omega_{\XX}$, which we then call the \emph{automorphic line bundle} on $\XX$.  We denote the total space of this line bundle by $\LL$ and the complement of its zero section by $\LL^\times$. The latter will play an important role in the various partial compactifications that we will discuss.

\subsection*{Domains of type I} The bounded symmetric domain  of type I$_{p,q}$ (in the Cartan classification) is defined by a finite dimensional complex vector space $W$ endowed with a nondegenerate Hermitian form $h:W\times W\to \CC$ of signature $(p,q)$. We take $G=\SU(W)$ and let $\XX$ be the open subset $\BB_W\subset\Gr_p(W)$ of positive definite subspaces of dimension $p$. The action of $\SU(W)$ on  $\BB_W$ is transitive and the stabilizer of $F\in \BB_W$ is contained in  the  product  $\U(F)\times \U(F^\perp))$ of compact unitary groups. This is indeed a maximal compact subgroup of $\SU(W)$. 

Given $F\in\BB_W$, then any other $F'\in\BB_W$ is the graph of a linear map 
$\phi: F\to F^\perp$. We identify $F^\perp$ with $V/F$ and observe that $\Hom (F, V/F)$ can be understood as the tangent space of the Grassmannian of $p$-planes in $W$ at the point defined by $F$. The condition that $h$ be positive on $F'$ translates into the following simple boundedness condition on  $\phi$: if 
$h_F:=h|F\times F$ and  $h_{F^\perp}:=-h|F^\perp\times F^\perp$ (so both are positive) then we require that $h_{F^\perp}(\phi (z),\phi (z))<h_F(z,z)$ for all $z\in F-\{ o\}$. This clearly implies that $\phi$ is bounded in 
$\Hom (F, F^\perp)$ and so this realizes $\BB_W$ as a bounded domain in this tangent space.

The tautological bundle of rank $p$ over $\Gr_p(W)$ restricts to one over $\BB_W$. We take as automorphic line bundle $\Lcal_W$ its $p$th exterior power. A straightforward check shows that this is a $(p+q)$th root of the canonical bundle of $\BB_W$.

The case $p=1$ is of special interest to us: then $\BB_W$ appears in $\Gr_1(W)=\PP(W)$ as a complex unit ball and the locus $W_+$ of $w\in W$ with $h(w,w)>0$ can be identified with the complement of the zero section $\LL^\times_W\subset\LL_W$. For $p=q=1$, we recover the usual upper half plane: there is a cobasis $(z_1,z_2)$ of $W$ such that 
$h(v,v')=-\sqrt{-1} (z_1(v)\overline z_2(v')- z_2(v)\overline z_1(v'))$ and  $\BB_W$ is then defined by $\im (z_1/z_2)>0$. In fact, the set defined by $h(v,v)=0$ defines a circle in $\PP(W)$ and there is a unique real form $V$ of $W$ for which this circle is the associated real projective line. The imaginary part of
$h$ defines a symplectic form $a$ on $V$ for which $\BB_W$ gets identified with $\HH_V$.

\subsection*{Domains of type II} A  domain of type II is given by a real vector space $V$ of finite even dimension $2g$ endowed with symmetric bilinear form $s: V\times V\to \RR$ of signature  $(g,g)$. Let $\DD_V$ 
be the set of $s_\CC$-isotropic $g$-dimensional subspaces $F\subset V_\CC$ that are $h$-positive, where 
$s_\CC$ resp.\ $h$  is the bilinear resp.\ hermitian extension of $s$. So this is naturally contained in a closed subset in a domain of type I$_{g,g}$: $\DD_V\subset \DD_{V_\CC}$, where $V_\CC$ stands for the complexification of $V$ endowed with the hermitian form $h$.  In fact, $\DD_V$ is obtained by intersecting $\DD_{V_\CC}$ with the $s$-isotropic locus in the Grassmannian $\Gr_g(V_\CC)$ (which is a nonsingular homogeneous subvariety).

The group $\Orth (V)$ acts transitively on $\DD_V$ and the stabilizer of any $F\in \DD_V$ restricts isomorphically to the unitary group $\U(F)$, a maximal compact subgroup of $\Orth (V)$. So $\DD_V$ is a symmetric space for the identity component of $\Orth (V)$.

\subsection*{Domains of type III} Perhaps the most familiar class of bounded symmetric domains are those  of type  III$_g$, usually presented as a tube domain and then known as a \emph{Siegel upper half space}. Such an object is naturally associated with a finite dimensional real vector space $V$ of dimension $2g$ endowed with a nondegenerate symplectic form $a: V\times V\to \RR$.  If $a_\CC: V_\CC\times V_\CC\to \CC$ denotes the complexification of the symplectic form, then the Hermitian form $h: V_\CC\times V_\CC\to \CC$ defined by  $h(v,v'):=\sqrt{-1}a_\CC (v, \overline v')$ has signature $(g,g)$. We take  $G=\Sp (V)$ and let  $\XX=\HH_V$ be the subset of the Grassmannian $\Gr_g(V_\CC)$ defined as follows: a $g$-dimensional subspace $F\subset V_\CC$ is in $\HH_V$ if $F$ is totally isotropic relative to $a_\CC$ and positive definite relative to $h$. The group $\Sp (V)$ acts transitively on $\HH_V$. The stabilizer of any $F\in \HH_V$ restricts isomorphically to the unitary group $U(F)$; it is a compact subgroup of $\Sp (V)$ and maximal for that property. That makes $\HH_V$ a symmetric space for  $\Sp (V)$. If $F\in \HH_V$, then $F$ defines a Hodge structure on $V$ of weight 1 polarized by $a$: take $V^{1,0}=F$ and $V^{0,1}=\overline F$. Thus $\HH_V$ parametrizes polarized Hodge structures on $V$ of this type. For this reason, the tautological bundle of rank $g$ over $\Gr_g(V_\CC)$ restricted to $\HH_V$ is often called the \emph{Hodge bundle}. We take for the automorphic  line bundle the sheaf $\Lcal_{V}$
of holomorphic sections of its $g$-th exterior power of the Hodge bundle. It is a $g$th root of the canonical bundle of $\HH_V$.

Any $F'\in \HH_V$ will not meet the $\overline F $ (for $h$ is negative on $\overline F$) and so must be the graph of a linear map $\phi :F\to \overline F$. Since $\overline F$  maps isomorphically to $V/F$ and $a$ identifies $ V/F$ with the complex dual of $F$, we can think of $\phi$ as a bilinear map $\Phi: (v,v')\in F\times F\mapsto a_\CC(v, \phi (v'))\in\CC$. The requirement that the graph of $\phi$ is totally isotropic relative to $a_\CC$ amounts to $\Phi$ being symmetric. In fact, the space of quadratic forms on $F$ may be regarded as the tangent space of $\HH_V$ at the point defined by $F$. The condition that  $h$ is positive on $F'$ amounts a boundedness condition on $\Phi$. Thus $\HH_V$ is realized as a bounded domain in the space of quadratic forms on $F$. For $g=1$, we thus recover the unit disk.

A type I domain can arise as a subdomain of type III domain as follows. With $(V,a)$ as above, suppose we are given a semisimple $\sigma \in \Sp (V)$ (in the sense that it decomposes $V_\CC$ into a direct sum of its 
eigenspaces: $V_\CC=\oplus_{\lambda} V^\lambda_\CC$). We determine what  its fixed point set in $\HH_V$, $\HH_V^\sigma$,  is like when nonempty.  Clearly, if $F\in \HH_V^\sigma$, then $\sigma\in U(F)\times U(\overline F)$ and so the  eigenvalues of $\sigma$ are 1 in absolute value and occur in complex conjugate pairs.  Write $F^\lambda:=F\cap V_\CC^\lambda$ and let $p_\lambda:=\dim F^\lambda$. 

First note that  $V^{\pm 1}_\CC$ is defined over $\RR$ and that $a$ is nondegenerate on the underlying real vector space $V^{\pm1}$. So $\HH_{V^{\pm1}}$ is defined and $F^{\pm 1}$ defines a point of $\HH_{V^{\pm1}}$.  Suppose $\lambda\not=\{ \pm1\}$. Then $a_\CC(v,v')=a_\CC (\sigma v, \sigma v')=\lambda^2a_\CC(v,v')$ and so $V^\lambda_\CC$ is $a_\CC$-isotropic. On the other hand, $a_\CC$ restricts to a nondegenerate pairing $V^\lambda_\CC\times  V^{\overline\lambda}_\CC\to \CC$. Since  $V^\lambda_\CC$ is the direct sum of $V^\lambda_\CC\cap F=F^\lambda$ and $V^\lambda_\CC\cap \overline F$ (which is just the complex conjugate of $F^{\overline \lambda}$), the hermitian form $h$  is nondegenerate on $V^\lambda_\CC$ and  of signature $(p_\lambda,p_{\bar\lambda})$ and $F^\lambda$ defines a point of the type I$_{p_\lambda,p_{\bar\lambda}}$  domain $\BB_{V^{\lambda}}$. Note that then $V^{\overline\lambda}_\CC=\overline V^\lambda_\CC$ has opposite signature $(p_{\bar\lambda}, p_\lambda)$ and that $F^{\overline\lambda}$ is the annihilator of $F^\lambda$ relative to $a_\CC$.  In particular,  $F^{\overline\lambda}$ is determined by  $F^\lambda$. We thus obtain an embedding 
\[
i: \HH_{V^{+1}}\times\HH_{V^{-1}}\times \prod_{\im(\lambda)>1} \BB_{V^\lambda_\CC}\hookrightarrow \HH_V.
\]
with image the fixed point set $\HH_V^\sigma$. This embedding is isometric up to rescaling   In particular, the image is a totally geodesic subdomain. The centralizer of $\sigma$ in $\Sp(V)$, $\Sp(V)_\sigma$, will act transitively on this product via an isomorphism $\Sp(V)_\sigma\cong\Sp (V^1)\times \Sp (V^{-1})\times \prod_{\im(\lambda)>0}\U (V^\lambda_\CC)$.
The pull-back of the automorphic line bundle $\LL_V$ can be equivariantly identified with the exterior tensor product of the automorphic line bundles of the factors.

Thus a polarized Hodge structure of type $(1,0)+(0,1)$ invariant under a semisimple symmetry gives rise to point in such a product decomposition. Our main interest will be in the factors that give rise to complex balls, i.e., those of type I$_{1,q}$.

Notice that we also have an inclusion of opposite type: a domain of type III$_g$ lies naturally in one of type I$_{g,g}$ as a closed orbit of the subgroup $\Sp (V)\subset \SU (V)$.

\subsection*{Domains of type IV}
This class is especially dear to algebraic geometers. A domain of type IV$_n$ is given by a real vector space $V$ of dimension $n+2\ge 2$ endowed with a symmetric bilinear form $s: V\times V\to\RR$ of signature $(2,n)$. Let $s_\CC: V_\CC\times V_\CC\to \CC$  denote the complexification  of $s$ and $h: V_\CC\times V_\CC\to \CC$ the Hermitian form defined by  $h(v,v')=s_\CC (v, \overline v')$. 

Denote by $V^+\subset V_\CC$ the set of $v\in V_\CC$ with $s_\CC(v,v)=0$ and $h(v,v)>0$. Its  projectivization $\PP(V^+)$ is an open subset in the nonsingular quadric in $\PP (V_\CC)$ defined by $s_\CC (v,v)=0$. If we assign to $[v]\in \PP(V^+)$ the oriented plane $P_{[v]}$ spanned by the real part and the imaginary part of $v$, then we obtain an identification of $\PP(V^+)$ with the Grassmannian  of positive definite oriented 2-planes in $V$. This makes it clear that  $\PP(V^+)$ has two connected components which are exchanged by complex conjugation. The orthogonal group $\Orth(V)$ acts transitively on $\PP(V^+)$ and the stabilizer of $[v]\in \PP(V^+)$ is $\SO (P_{[v]})\times \Orth (P_{[v]}^\perp)$. Its intersection  with the identity component $\Orth(V)^\circ$ is $\SO (P_{[v]})\times \SO (P_{[v]}^\perp)$ and this is a maximal compact subgroup of $\Orth(V)^\circ$. So each connected component of $\PP(V^+)$ is a bounded symmetric domain for $\Orth(V)^\circ$. We now let $\LL^\times$ be a connected  component of $V^+$ and write $\DD_V$ for its projectivization $\PP(\LL^\times)$. From the preceding it is clear that 
$\DD_V$ is a bounded domain. The \emph{automorphic line bundle} $\Lcal_V$ on $\DD_V$ is simply the restriction of the tautological bundle over $\PP (V_\CC)$ so that $\LL^\times$ is indeed its total space minus its zero section. It is an $n$th root of the canonical bundle of $\DD_V$. We denote by $\Orth^+(V)$ the subgroup of 
$\Orth(V)$ which leaves $V^+$ invariant.

Any $[F]\in\DD_V$ determines a bounded realization of  $\DD_{V}$  in the usual  manner: any $[F']\in\DD_{V}$ is the graph of a linear map $\phi=(\phi', \phi''): F\to (F\cap F^\perp)\oplus \overline F$. The condition that this graph is $s_\CC$-isotropic  implies that $\phi''$ is determined by $\phi'$,   for we then must have $s_\CC(v, \phi''(v))=-s_\CC(v,v)-s_\CC(\phi' (v),\phi' (v))$. If we identify $F\cap F^\perp$ with $F^\perp/F$, then we see that $\phi$ is given by an element of $\Hom (F, F^\perp/F)$, which is just the tangent space of $\DD_V$ at $[F]$ and thus  $\DD_{V}$ embeds in this tangent space.  The  condition that $h$ be positive on $F'$ implies the boundedness of the image.

We have already seen the cases $n=1,2$ in a different guise: if $V_1$ is of real dimension two and endowed with a real symplectic form $a: V_1\times V_1\to\RR$, then a symmetric bilinear form $s$ on $V_1\times V_1$ is defined by $s(v\otimes w, v'\otimes w')=a(v,v')a(w, w')$. Its signature is $(2,2)$, and  $(F, F')\mapsto F\otimes F'$ 
defines an isomorphism of $\HH_{V_1}^2$ onto a connected component of $\DD_{V_1\otimes V_1}$ that is compatible with an isogeny of $\Sp (V_1)\times\Sp (V_1)$ onto $\Orth (V_1\otimes V_1)^\circ$. Here $\Lcal_{V_1\otimes V_1}$ pulls back to the exterior product $\Lcal_{V_1}\boxtimes\Lcal_{V_1}$. If we restrict this map
to the diagonal of $\HH_{V_1}^2$, then we land in the intersection of $\DD_{V_1\otimes V_1}$ with the subspace 
$\PP(\sym^2V_1)\subset \PP(V_1\otimes V_1)$. The restriction of $s$ to  $\sym^2V_1$ has signature $(1, 2)$ and
we thus get an identification of $\HH_{V_1}$ with a component of $\DD_{\sym^2V_1}$. This is compatible with an isogeny of $\Sp (V_1)$ onto $\Orth (\sym^2V_1)^\circ$. Notice that  $\Lcal_{\sym^2V_1}$ pulls back to $\Lcal^{\otimes 2}_{V_1}$. 

Any $F\in \DD_V$ gives rise to a decomposition $V_\CC=F\oplus (F\oplus \overline F)^\perp\oplus \overline F$, which we can think of as defining a Hodge structure of weight 2 polarized by $s$ with $F=V^{2,0}$, $V^{1,1}=(F\oplus \overline F)^\perp$ and $\overline F=V^{0,2}$.
\\

In this context  ball domains may arise in a similar manner. Suppose $\sigma\in \Orth^+(V)$ is semisimple and $\DD_V^\sigma$ is nonempty. Clearly, if $[v]\in\DD_V^\sigma$, then  $v\in V_\CC^\lambda$ for some eigenvalue $\lambda$ of $\sigma$. 
Since $\sigma\in \SO (P_{[v]})\times \SO (P_{[v]}^\perp)$, all the eigenvalues of $\sigma$ lie on the complex unit circle with the ones $\not=\pm 1$ appearing in complex conjugate pairs.
When $\lambda\in \{\pm 1\}$, $V_\CC^\lambda$ is defined over $\RR$ and must have signature $(2, n')$ for some integer $n'\ge 0$. This yields a totally geodesic embedding  $\DD_{V^{\lambda}}\cong \DD_V^\sigma\subset \DD_V$ and $\Orth(V)_\sigma$ acts via a homomorphism $\Orth(V)_\sigma\to \Orth(V^{\lambda})$. This homomorphism is onto and has compact kernel. 

For $\im (\lambda) >0$,  $h$ will have the same signature on
$V_\CC^{\overline\lambda}=\overline V_\CC^{\lambda}$ as on $V_\CC^\lambda$. It follows that this signature must be of the form $(1,n')$ for some integer $n'\ge 0$ and that all other eigenspaces are negative definite. Hence $\DD_V^\sigma$ is connected and is embedded as a totally geodesic submanifold. The group $\Orth^+(V)_\sigma$ acts via a homomorphism $\Orth^+(V)_\sigma\to \U(V_\CC^\lambda)$ that is onto and has compact kernel.  

\subsection*{The exceptional domains} Two isomorphism types remain: the domains  that are the symmetric spaces of real forms of simple Lie groups of type $E_6$ and $E_7$. They have dimension 16 and 27 respectively. The bounded domain of type $E_7$ supports a variation of Hodge structure of weight 3 and rank 56 as follows.  Let $V$  be a real vector space of dimension $56$. There exists in $\PP (V_\CC)$ a  quartic hypersurface $Q$ defined over $\RR$ with the property that its singular locus $S\subset Q$  is nonsingular  of dimension $27$ and is homogeneous for the subgroup $\Gcal$ of $\SL(V_\CC)$ which preserves $Q$. This group $\Gcal$ is then necessarily of type $E_7$ and defined over $\RR$, $V_\CC$ is an irreducible representation of it  and  $S$  the highest weight orbit. The 54-dimensional hypersurface $Q$ can be recovered from $S$ as the union of the projectively completed tangent spaces of $S$.  Moreover, there is a nondegenerate  antisymmetric form $a: V\times V\to \RR$ with the property that the  projectively completed tangents of $S$ are projectivized Lagrangian subspaces of $a_\CC$. The form $a$ is unique up scalar and hence $\Gcal$-invariant.
We thus have defined for  every $s\in S$ a flag $\{ 0\}=F^{4}_s\subset F^{3}_s\subset F^{2}_s\subset F^{1}_s\subset F^{0}_s=V_\CC$ as follows: $F^{3}_s\subset V_\CC$  is the one dimensional subspace associated with $s$ and $F^{2}_s\supset F^{3}_s$ the  subspace of dimension 28 defined by the tangent space $T_sS$.  This is Lagrangian subspace for $a_\CC$ and so the  annihilator $F_s^1$ of $F^{3}_s$  relative to $a_\CC$ will contain  $F^{2}_s$.
The stabilizer $\Gcal_s$ acts faithfully on $F^{2}_s$: the action on $F^3_s\times F^2_s/F^3_s$ is via its Levi quotient, a group isomorphic to $\CC^\times$ times a group of type $E_6$ and the unipotent radical of $\Gcal_s$ is a vector group that can be identified with $\Hom (F^2_s/F^3_s,F^3_s)$.  

The identity component of the real part of $\Gcal$, $G:=\Gcal (\RR)^\circ$, has an open orbit $\XX$  in $S$  whose points $o\in \XX$ have the property that the hermitian form  associated with $a$ is positive on $F^{3}_o$ and of signature $(1, 27)$ on $F^2_o$. This  defines a variation of polarized Hodge structure of weight $3$ over $\XX$ with Hodge numbers $h^{3,0}=h^{0,3}=1$ and $h^{2,1}=h^{1,2}=27$.   It is effective  in the sense that $T_o\XX$ maps isomorphically to $\Hom (F^{3,0}_o, F^{2,1}_o/F^{3,0}_o)$. But this variation of  Hodge structure has not yet been encountered in algebraic geometry.

It seems likely that we can obtain a bounded symmetric domain of type $E_6$ in $\XX$ as a locus with additional symmetry. For instance if we have a $s\in S$ such that $h$ is nonzero on $F^{3}_s$  and of signature $(11, 16)$ in
$F^2_s/F^3_s$, then $G_s$ acts on $H^{2,1}_s:=F^2_s\cap (F^3_s)^\perp$ through a real group of  type $E_6$  and I suspect that for the  correct sign of $h$ on $F^3_s$,  $\XX\cap \PP(F^2_s)$ is its symmetric domain. 



\subsection*{Hermitian domains admitting complex  reflections}The complex balls and the domains of type IV are the only bounded symmetric domains which admit totally geodesic subdomains of codimension one (or equivalently, admit complex reflections). In one direction this is easy to see: if $(W,h)$ is a Hermitian form of hyperbolic signature $(1,n)$, then any complex hyperplane section of $\BB_W$ has this property: if $W'\subset W$ is a subspace of complex codimension one such that $\BB_W\cap \PP (W')$ is nonempty, then $W'$ will have signature $(1, n-1)$ and $\BB_{W'}=\BB_W\cap \PP (W')$ is the fixed point set of a copy of  $\U (1)$ (identify $\U(1)$ with $\U ((W')^\perp)$). In particular, $\BB_{W'}$ is geodesically embedded. 

A similar argument shows that if $(V,s)$ defines a domain of type IV as above, then any complex hyperplane section of $\DD_V$ defined over the reals has the desired property: if $V'\subset V$ is of real codimension one and of signature $(2, n-1)$, then $V'$ is the fixed point set of a reflection  $\sigma\in \Orth(W)$ and so $\DD_{V'}=\DD_V\cap \PP (V'_\CC)=\DD_V^\sigma$. 

\section{The Baily-Borel package}\label{section:bbpackage}

\subsection*{Arithmetic structures} The notion of an arithmetic subgroup of an algebraic group requires the latter to be defined over a number field. A trick, the so-called restriction of scalars, shows that here is in fact little or no loss in generality in assuming that this number field is $\QQ$. In the cases $\Sp (V)$ and $\Orth (V)$ this is simply accomplished by assuming that $(V,a)$ resp.\ $(V,s)$ is defined over $\QQ$. An arithmetic subgroup $\Gamma$ of $\Sp (V)$ resp.\ $\Orth (V)$ is then a subgroup which fixes a lattice $V_\ZZ\subset V_\QQ$ and is of finite index in $\Sp (V_\ZZ)$ resp.\ $\Orth (V_\ZZ)$. Such a group is discrete in the ambient  algebraic group and hence acts properly discontinuously on the associated Hermitian domain. 

Examples for  the case of complex ball can be obtained as in the discussion above by taking a  semisimple $\sigma\in \Sp (V_\ZZ)$ resp.\  $\sigma\in \Orth^+(V_\ZZ)$ that has fixed points in the associated Hermitian domain. Since $\sigma$ lies in a unitary group and fixes a lattice, it must have finite order, $m$ say so that its eigenvalues  are $m$th roots of unity. 
If $\mu_m^\pr$ denotes the set of primitive roots of unity, then $\oplus_{\lambda\in \mu_m^\pr} V^{\lambda}_\CC$ is defined over $\QQ$. So the subgroup $G$ of $\Sp (V)_\sigma$ resp.\ $ \Orth (V_\ZZ)$ which acts as the identity on the remaining eigenspaces is defined over $\QQ$ and the  subgroup $G(\ZZ)$ which stabilizes the lattice $V_\ZZ$ is arithmetic in $G$. 
Suppose first we are in the symplectic case. If for all $\lambda\in \mu_m^\pr$, the form $h$ is definite on $V_\CC^\lambda$ except for a complex conjugate pair  $(\lambda_o, \overline\lambda_o)$ with $h$ of signature $(1,n)$ on $V_\CC^{\lambda_o}$ (and hence of signature $(n,1)$ on $V_\CC^{\overline\lambda_o}$), then $G$ acts properly on $\BB_{V_\CC^{\lambda_o}}$ and hence $G(\ZZ)$ acts on this domain properly discretely. A similar argument applies to the orthogonal case. The situation simplifies considerably when  $m=3,4,6$, for then there are only  two primitive $m$th roots of unity (namely $\zeta_m:=\exp(2\pi\sqrt{-1}/m)$ and its conjugate $\overline \zeta_m:=\exp(-2\pi\sqrt{-1}/m)$). So the above decomposition reduces to $V^{\zeta_m}_\CC\oplus V^{\overline\zeta_m}_\CC$ and $G$ has no compact factors unless $G$ itself is compact.

The more intrinsic approach is to interpret an arithmetic structure on $(W,h)$ as given by a lift $(W_K,h_K)$ over a CM-field $K/\QQ$ with the property that for every embedding $\iota: K\to \CC$ the corresponding Hermitian form on $W_K\otimes_\iota \CC$ is definite except for a complex conjugate pair, for which the signatures  are $(n,1)$ and $(1,n)$.

\subsection*{The Baily-Borel package for locally symmetric varieties}
Let $(G, \LL)$ be one of the examples above. To be precise: $G$ is an algebraic group defined over $\RR$ and $G(\RR)$ acts properly and transitively on a  complex manifold $\XX$ that is a symmetric domain for $G(\RR)$ and $\LL$ is the total space of a $G$-equivariant line bundle over $\XX$ that is an equivariant positive root of the canonical bundle with $G(\RR)$-action. We denote by $\LL^\times\subset\LL$ the complement of the zero section.

Now suppose $G$ is defined over $\QQ$ and we are given an arithmetic subgroup  $\Gamma\subset G(\QQ)$. Then
$\Gamma$ acts properly discretely on $\XX$ (and hence on $\LL^\times$). The action of $\Gamma$ on $\LL^\times$ commutes with the obvious action of $\CC^\times$ so that we have an action of $\G\times\CC^\times$ on $\LL^\times$. 
The Baily-Borel package that we are going to state below can be entirely phrased in terms of this action. It begins with the observation that  we can form the $\G$-orbit space of the $\CC^\times$-bundle $\LL^\times \to \XX$ in the complex-analytic category of orbifolds to produce an $\CC^\times$-bundle $\LL^\times_\G\to\XX_\G$. For $d\in\ZZ$, we denote by $A_d(\LL^\times)$ the linear space of holomorphic functions $f:\LL^\times\to \CC$ that are homogeneous of degree $-d$ (in the sense that $f(\lambda z)=\lambda^{-d}f(z)$) and subject to a growth condition (which in many cases is empty). This growth condition is such that $A_\pt (\LL^\times)$ is closed under multiplication (making it a graded $\CC$-algebra) and invariant under $\G$.
\begin{description}
\item[Finiteness] The graded $\CC$-subalgebra of $\G$-invariants, $A_\pt (\LL^\times)^\G$ is finitely generated with generators in positive degree. So we have defined the normal weighted homogeneous (affine) cone $(\LL^\times_\G)^\bb:=\spec A_\pt (\LL^\times)^\G$ whose base is the normal projective variety $\XX_\G^\bb:=\proj (A_\pt (\LL^\times)^\G)$. 
\item[Separation] The elements of $A_\pt (\LL^\times)^\G$ separate the $\G$-orbits in $\LL^\times$ so that the natural maps $\LL^\times_\G\to (\LL^\times_\G)^\bb$ and $\XX_\G\to \XX_\G^\bb$ are injective.
\item[Topology] The underlying analytic space of $(\LL^\times_\G)^\bb$ is obtained as $\G$-orbit space of a ringed space $((\LL^\times)^\bb, \Ocal_{(\LL^\times)^\bb})$ endowed with a $\G\times \CC^\times$-action and 
a $\G\times \CC^\times$-invariant partition into complex manifolds such that $\Ocal_{(\LL^\times)^\bb}$ is the sheaf of continuous complex valued functions that are holomorphic on each part. This partition descends to one of  $(\LL^\times_\G)^\bb$ into locally closed subvarieties which has $\LL^\times_\G$ as an open-dense member. (Thus $\XX_\G^\bb$ comes  with stratification into locally closed subvarieties which has the image of $\XX_\G$ as an open-dense member. In particular,  $\LL^\times_\G$ resp.\ $\XX_\G$
acquires the structure of a quasi-affine resp.\ quasi-projective variety.)
\end{description}
In fact, $(\LL^\times)^{\bb}$ and $A_\pt (\LL^\times)$ only depend  on $G(\QQ)$ and the $\G\times \CC^\times$-action is the restriction of one of  $G(\QQ)\times \CC^\times$, but in generalizations that we shall consider the latter group is no longer acting.
We shall make the topological part somewhat explicit in some of the classical cases discussed above.

\subsection*{Ball quotients} These are in a sense the simplest because $\BB_\G^\bb$ adds to $\BB_\G$ a finite set. Denote by $\Ical$ the collection of isotropic subspaces of $W$ that is left pointwise fixed by an element of $G(\QQ)$ of infinite order (one can show that we can  then take this element to lie in $\G$). Note that every $I\in \Ical$ is of dimension $\le 1$. It is known that $\G$ has finitely many orbits in $\Ical$.

For every $I\in\Ical$, we have a projection $\pi_{I^\perp} : \LL^\times\subset W\to W/I^\perp$. When $I$ is of dimension one, the image of $\pi_{I^\perp}$  is the complement of the origin. As a set 
$(\LL^\times)^\bb$, is the disjoint union of $\LL^\times$ and its projections $\pi_{I^\perp}(\LL^\times)$. It comes with an evident 
action of $G(\QQ)\times\CC^\times$ and the topology we put on it makes this action topological. It has a conical structure with vertex the stratum defined by  $I=\{ 0\}$. So $(\LL^\times)^\bb$ adds to $\LL^\times$ a finite number of $\CC^\times$-orbits, one of which is a singleton. In particular, $\BB^\bb_\G-\BB_\G$ is finite.

\subsection*{Quotients of orthogonal type} Here $(V,s)$ is a real vector space endowed with a nondegenerate quadratic form   defined over $\QQ$ and of signature $(2,n)$. Denote by $\Ical$ the collection of isotropic subspaces of $V$ defined over $\QQ$. 

For every $I\in\Ical$, we have a projection $\pi_{I^\perp} : \LL^\times\subset V_\CC\to V_\CC/I^\perp_\CC$. When $\dim I=2$, then the image consists of elements in $V_\CC/I^\perp_\CC$ that under the identification $V_\CC/I^\perp_\CC\cong\Hom_\RR (I, \CC)$ have real rank 2. Since $\LL^\times$ is connected, they all define the same orientation of $I$ and indeed, $\pi_{I^\perp}$ is the space $\iso^+(I, \CC)$ of orientation preserving real isomorphisms. The $\CC^\times$-orbit space of $\iso^+(I, \CC)$ is a copy of the upper half plane and parametrizes  all the complex structures on $I$ compatible with the orientation of $I$. The $\G$-stabilizer $\G_I$ of $I$ acts on 
$\iso^+(I, \CC)$ via an arithmetic subgroup of $\SL (I)$  and the orbit space of  $\PP(\iso^+(I, \CC))$ is  a modular curve. These modular curves are noncompact, but their missing points are supplied by the $I\in\Ical$ of dimension one.
When $\dim I=1$ the image is the punctured line $(V/I^\perp)_\CC-\{ 0\}\cong \Hom_\RR(I, \CC)-\{ 0\}$ and for $I=\{ 0\}$ we get a singleton. 

Now  $(\LL^\times)^\bb$ is the disjoint union of $\LL^\times$ and its projections  $\pi_{I^\perp} (\LL^\times)$ and endowed with a 
$\Orth^+(V_\QQ)\times\CC^\times$-equivariant topology. So $\DD^\bb_\G-\DD_\G$ consists of a union of finitely many  modular curves and finitely many points. We will describe this construction in more detail below.

\subsection*{Moduli of polarized abelian varieties} 
Assume given a symplectic vector space  $(V,a)$ defined over $\QQ$. Since two symplectic forms in the same number of variables are equivalent, there is essentially one such space when its dimension is given. Recall that here the symmetric domain is $\HH_V$, the set of $g$-dimensional subspaces $F\subset V_\CC$ that are totally isotropic relative to $a_\CC$ and positive definite relative to $h$. Denote by $\Ical$ the collection of isotropic subspaces of $V$ defined over $\QQ$. Given $I\in \Ical$, then $I^\perp\supset I$ and  $I^\perp/I$ is a symplectic vector space over $\QQ$. It is then easily seen that $\Sp (V_\QQ)$ acts transitively on the collection of $I\in \Ical$ of prescribed dimension. Given $I\in \Ical$, then for every $F \in\HH_V$,   $F\cap I_\CC=\{ 0\}$, $F\subset V_\CC\to V_\CC/I^\perp_\CC\cong\Hom_\RR(I, \CC)$ is onto and $F\cap  I^\perp_\CC$ maps isomorphically onto an element $F_I$of $I^\perp/I$. We have thus defined
a projection $\HH_V\to\HH_{I^\perp/I}$. This projection  lifts to the automorphic line bundle as follows: we observe that $\det (F)\cong \det(F\cap  I^\perp)\otimes_\RR\det (V/I^\perp)\cong  \det(F_I)\otimes_\RR \det (I)^{-1}$ and this suggests that the natural line bundle of $\HH_{I^\perp/I}$ must be twisted by $\det (I)$.

\subsection*{Baily-Borel theory and its variations in a nutshell} We begin with sketching the nature of the Baily-Borel extension $\XX^\bb$ for $G$ as above (so  defined over $\QQ$) in general terms. We assume for simplicity that
$\Gcal$ is almost-simple as a $\QQ$-group. A central role is played by the collection $\Pcal_{\max}(G(\QQ))$ of maximal parabolic subgroups of $G$ defined over $\QQ$. Indeed, the  extension in question requires a thorough understanding of the structure of such subgroups and so we address this first.

Let  $P\subset G$ be a maximal (for now, only real) parabolic subgroup of $G$. Then its unipotent  radical $R_u(P)\subset P$ is rather simple: if $U_P\subset 
R_u(P)$ denotes its center (a vector group), then $V_P:=R_u(P)/U_P$ is also a (possibly trivial) vector group. It is now best to pass to the associated Lie algebras (where we shall follow the custom of denoting these in the corresponding Fraktur font). The Lie bracket then defines an antisymmetric bilinear map $\vfrak_P\times \vfrak_P\to \ufrak_P$. This map is of course equivariant with respect to the adjoint  action of $P$
on these vector spaces. One finds that $\ufrak_P$ contains a distinguished  open orbit $C_P$ of $P$ with remarkable properties: it is a strictly convex cone and if we exponentiate  $\sqrt{-1}C_P$ to a semigroup in $\Gcal (\CC)$, then this semigroup leaves $\XX$ invariant (think of the upper half plane in $\CC$ that is invariant under the  semigroup of translations $\sqrt{-1}\RR_{>0}$). Moreover, for  every linear form $\ell$ on $\ufrak_P$ which is positive on $\overline{C}-\{ 0\}$, the composite of
$\vfrak_P\times \vfrak_P\to \ufrak_P\xrightarrow{\ell}\RR$ is nondegenerate so that dividing $R_u(\pfrak)$ out by the kernel of $\ell$ yields a Heisenberg algebra (which is trivial when $\vfrak_P$ is). The Lie algebra $\ufrak_P$ (or $C_P$ for that matter) completely determines $P$, for we can recover $P$ as the $G$-stabilizer of  $\ufrak_P$ (resp.\ $C_P$). In fact, this gives rise to an interesting partial order on the collection $\Pcal_{\max}(G(\RR))$ of maximal parabolic subgroups of $G$:  we define $Q\le P$ in case $\ufrak_Q\subset\ufrak_P$ (in $\gfrak$) and this last property is equivalent to $C_Q\subset \overline C_P$. The various $C_Q$'s are pairwise disjoint and one finds that the 
union of the $C_Q$ over all $Q$ with $Q\le P$ makes up all of $\overline C-\{ 0\}$.

The action of $P$ on $\ufrak_P\oplus\vfrak_P$ is through its Levi quotient  $L_P:=P/R_u(P)$, the latter acting  with finite kernel.  We write $M^\hor_P\subset L_P$ for the kernel of the action of $L_P$ on $\ufrak_P$. This is a semi-simple subgroup. In fact, it is an almost  direct factor of $L_P$, for the centralizer of $M^\hor_P$ in  $L_P$ is a reductive group which supplements  $M^\hor_P$ in $M_P$ up to a finite group. This reductive group, which we denote by  $L^\ver_P$ decomposes naturally as  $A_P.M^\ver_P$ with $A_P$ isomorphic to $\RR^\times$ and $M^\ver_P$ semisimple. (So the preimage of $A_P$ in $P$ is the radical $R(P)$ of $P$ and $P/R(P)$ is the almost product  $M^\hor_P.M^\ver_P$. The group $A_P\cong \RR^\times$ acts by a nontrivial character---in fact, by squares---on the cone $C_P$ and the real projectivization of $C_P$ is a symmetric space for  $M^\ver_P$.)

It turns out that the symmetric space of $M^\hor_P$ is a symmetric domain and that this domain is naturally obtained as a holomorphic quotient $\XX(P)$ of $\XX$.  In fact, $\XX\to \XX(P)$ is defined by the property that its fibers are the maximal orbits in $\XX$ of the semi-group $R_u(P)+\exp (\sqrt C_P)$ in $\Gcal (\CC)$. This is essentially  (but put in more abstract form) what is known as a realization of $\XX$ as as a Siegel domain. (The domain $\XX(P)$ can  naturally be realized in the boundary of $\XX$ with respect to the Harish-Chandra embedding and this is compatible with the opposite partial order above: $\XX(P)$  lies in the closure of $\XX(Q)$ precisely when $Q\le P$ and we thus obtain the decomposition of the boundary of  $\XX$ into what are called its boundary components.  But the topology that comes from this embedding  is not the one that matters here.) 

We now assume that $P$ is defined over $\QQ$. This of course will put a $\QQ$-structure on all the associated maps and spaces we encountered above. In order to describe a  topology of the disjoint union of $\XX$ and $\XX(P)$, it is convenient to fix  an arithmetic subgroup $\G\subset \G(\QQ)$, so that $\G_P:=P\cap \G$ is 
 an arithmetic subgroup of $P$. In particular, $\G\cap R_u(P)$ is  an extension of a lattice by a lattice. The image of $\G_P$ in $\GL (\ufrak_P)$ is a discrete subgroup (denoted $\G (P)$) which preserves $C_P$.  This is a group  which acts properly discretely  on $C_P$. In fact, if we let $C_{P}^+$ be the union of  all the $C_Q$ with $Q\le P$ with $Q\in\Pcal_{\max}(G(\QQ))$ and the origin $\{0\}$, then $C_{P}^+$ is the cone spanned by the set of  rational vectors in $\overline C_P$ and  $\G (P)$ has in $C_P^+$ a fundamental domain that is a rational polyhedral cone (i.e., the convex cone spanned by a finite number rational vectors in $C_P^+$). 
 
 The topology we impose on $\XX\sqcup \XX(P)$ will be such that every subset $\Omega\subset \XX$ that is invariant under both the semigroup $\sqrt{-1}C_P$ and the group  $R_u(P).\G (L^\ver_P)$ will have its image $\pi_P(\Omega)$ in $\XX(P)$ in its closure; here $\G (L^\ver_P)$ stands for the  intersection of the image of $\G_P\to L_P$ with $L^\ver_P$. So the topology we choose is generated by the open subsets of  $\XX$  and the subsets of the form $\Omega\sqcup \pi_P(\Omega)$, where $\Omega$ runs over the open subsets of $\XX$ with the invariance properties mentioned above. This is independent of our choice\footnote{Since we get the same topology if we replace $R_u(P).(\G_P\cap L^\ver_P)$ by $R_u(P).(\G_P\cap M^\ver_P)$, one may be tempted to think that no choice was needed, and  replace $R_u(P).(\G_P\cap L^\ver_P)$ by  $R_u(P).M^\ver_P$. But this topology  will in general differ from the topology that we want and which works.}  of $\G$. The topology induced on $\XX(P)$ is easily seen to be the given topology and to be such that we can extend $\pi_P$ to a continuous retraction of $\XX\sqcup \XX(P)$ onto  $\XX(P)$. 
 
We now put things together as follows. Let
\[
\XX^\bb:=\XX\sqcup \bigsqcup_{P\in\Pcal_{\max}(G(\QQ))}\XX(P).
\]
We equip this disjoint union with the topology generated by the open subsets of $\XX$ and the subsets $\Omega^{\bb_P}$ of the following type: for any $P\in\Pcal_{\max}(G(\QQ))$ and open subset $\Omega\subset \XX$ invariant under both the semigroup $\sqrt{-1}C_P$ and the group  $R_u(P).(\G\cap M^\ver_P)$, we let
\[
\Omega^{\bb_P}:=\Omega \sqcup \bigsqcup_{Q\in\Pcal_{\max}(G(\QQ)); Q\le P} \pi_Q(\Omega).
\]
This strangely defined space and its $\G$-orbit space is best studied at each of its `corners' separately. By this we mean that we fix a $P\in \Pcal_{\max}(G(\QQ))$ and restrict our attention to the star of $\XX(P)$, $\XX^{\bb_P}$.
Then the collection  $\{\pi_Q\}_{ Q\le P}$ defines a continuous retraction $\XX^{\bb_P}\to \XX(P)$ that is equivariant with respect to the 
$P(\QQ)$-action.  
An important and useful feature of this topology is that every $x\in \XX(P)$ has a basis of $\G_x$-invariant neighborhoods in $\XX^\bb$. Now the group $\G_P$ acts on $\XX(P)$ via a quotient $\G (\XX(P))$ that acts properly and as an arithmetic group on $\XX(P)$.  In particular, $\G_x$ contains the group $\G_P^\ver:=\G\cap R_u(P)M^\ver_P$ (which acts as the identity on $\XX(P)$) as a subgroup of finite index. So the main part about  understanding the $\G$-orbit space of $\XX^\bb$ near the image of $x$ is in understanding  the $\G_P^\ver$-orbit space of $\XX^{\bb_P}$: the latter must be shown to produce a normal analytic space with a proper $\G_P/\G_P^\ver$-action which comes with an equivariant analytic contraction onto $\XX(P)$. The natural map 
$\XX^{\bb_P}_{\G_P}\to \XX^\bb_\G$  may then be understood as providing a chart of $\XX^\bb_\G$ near the image of $\XX(P)$. But this is where we stop and refer to the examples discussed in this paper or to the original paper \cite{bb} to get an idea how to understand $\XX^{\bb_P}$.

In this paper, we are also concerned with the `automorphic cone' over the Baily-Borel compactification. This is an affine cone which contains $\LL^\times_\G$ as an open-dense subvariety.  It is defined in much the same way as $\XX^\bb_\G$, namely as the $\G$-orbit space of a Baily-Borel extension $(\LL^\times)^\bb$  of  $\LL^\times$, for  its boundary decomposes as a disjoint union of strata $\LL^\times(P)$ (defined as the orbit space of the action of the semigroup $R_u(P)+\exp (\sqrt C_P)$)  plus a vertex.
\\

\subsubsection*{Mumford's toroidal modification in a nutshell.} At this point, it is however relatively easy to say what data are needed for Mumford's toroidial compactification as described in \cite{amrt} and how this compares with the Baily-Borel construction. Recall that the union of the $C_P$'s in $\gfrak$ is a disjoint one and that the  closure of a member $C_P$ in $C(\gfrak(\QQ))$ is 
$C_P^+$.  We denote that union by $C(\gfrak(\QQ))$. The extra ingredient  we then need  is a $\G$-invariant refinement  $\Sigma$ of this decomposition of  $C(\gfrak(\QQ))$ into relatively open rational convex polyhedral cones (convex cones spanned by a finite subset of $\gfrak (\QQ)$). In a sense this greatly simplifies the  discussion, because the role of $Z_\G (\XX(P))$ is now taken over by the much smaller and simpler lattice $\G\cap U_P$ and the hard-to-understand $\XX^\bb_{Z_\G (\XX(P))}$ is replaced by a much easier understand  toroidal  extension of $\XX_{U_L}$, the torus in question being 
$T_P:=U_P(\CC)/\G\cap U_P$: we first enlarge $\XX$ to a complex manifold that we denote $U_P(\CC).\XX$, because it comes with transitive action of the subgroup $P.U_P(\CC)$ of $\Gcal (\CC)$ such that $\XX$ is an open $P$-orbit in $U_P(\CC).\XX$. Then $\XX_{\G\cap U_P}$ is open in   $(U_P(\CC).\XX)_{\G\cap U_P}$ and the latter is a principal $T_P$-bundle. The decomposition of $\Sigma| C_P^+$ defines a relative torus embedding  of $(U_P(\CC).\XX)_{\G\cap U_P}$. The interior of the closure of $\XX_{\G\cap U_P}$ in this relative  torus embedding is the  toric extension in question. It is an analytic variety with toric singularities on which the group $\G_P/(\G_P\cap U_P)$ acts properly discontinuously. 

From our point of view  it is more natural (and also closer the Satake-Baily-Borel spirit) to do this construction before dividing out by $\G_P\cap U_P$, that is, to define a $\G$-equivariant  extension $\XX^\Sigma$ of $\XX$. Although this brings us, strictly speaking, again outside the setting of analytic spaces, this is still simple enough. For every $\sigma\in\Sigma$ we can form a holomorphic  quotient  
$\pi_\sigma: \XX\to \XX(\sigma)$ whose fibers are maximal orbits for the semigroup $\exp (\la \sigma\ra_\RR+\sqrt{-1}\sigma)$. We let $\XX^\Sigma$  be the disjoint union of  $\XX$ and the $\XX(\sigma)$ and equip this union with the topology generated by the open subsets of $\XX$ and those of the form $\Omega^{\bb, \sigma}$: here 
$\sigma\in\Sigma$,  $\Omega\subset \XX$ is an open subset  invariant under the semigroup $\exp (\la \sigma\ra_\RR+\sqrt{-1}\sigma)$, and
\[
\Omega^{\bb, \sigma}:=\bigsqcup_{\tau\in\Sigma; \tau\subset \overline\sigma} \pi_\tau(\Omega)
\]
(note that $\Omega$ appears in this union for $\tau=\{0\}$).   When $\sigma\subset C_P$, then $\XX(\sigma)$ fibers naturally over $\XX(P)$. These projections combine together to define a continuous $\G$-equivariant  map $\XX^\Sigma\to \XX^\bb$, such that  the resulting map $\XX^\Sigma_\G\to \XX^\bb_\G$ is a morphism in the analytic category.
\\

\subsubsection*{Between the Baily-Borel and the Mumford constructions.}
A variation  discussed in this paper interpolates between these two. In the cases at hand\footnote{In general we would also have to specify for every $\sigma\in\Sigma |C_P$ a subspace $\vfrak_\sigma\subset\vfrak_P$ that is subject to certain properties, but such a choice is automatic in
the case  the case that we shall consider here (namely that of an `arithmetic arrangement').}
 we accomplish this by simply relaxing the condition on $\Sigma$: we just ask that it is a $\G$-invariant  decomposition of $C(\gfrak(\QQ))$ into relatively open \emph{locally} rational convex polyhedral cones (i.e., convex cones with the property that if we intersect them with rational convex polyhedral cone  we get a rational convex polyhedral cone; see \cite{looijcone} for details). Each $C_P$ is a locally  rational convex polyhedral cone and indeed, the coarsest choice is the decomposition of $C(\gfrak(\QQ))$ into $C_P$'s and will yield the Baily-Borel extension.

\subsection*{Arithmetic arrangements} As we noted,  a domain of type I$_{1,n}$ or IV   has totally geodesic complex hypersurfaces, which appear in our description as hyperplane sections. For such a domain $\XX$, and where $G$ is defined over $\QQ$, we make the following definition. 

\begin{definition}\label{def:}
An \emph{arithmetic arrangement} in $\XX$ is a collection $\Hcal$ of totally geodesic hypersurfaces that are also hermitian symmetric with the property that the $G$-stabilizer of each member is defined over $\QQ$ and that there is an arithmetic subgroup $\G\subset G(\QQ)$ such that $\Hcal$ is a finite union of $\G$-orbits.
\end{definition}

Then the collection $\{H\}_{H\in\Hcal}$ is locally finite so that their union $\Delta_\Hcal:=\cup_{H\in\Hcal}H$ is closed in $\XX$.  For $\G\subset G(\QQ)$ as above,  the image of this union in the $\G$-orbit space 
$\XX_\G$ is a hypersurface $(\Delta_{\Hcal})_\G\subset  \XX_\G$. This hypersurface is locally given by a single equation (it is $\QQ$-Cartier).  We call such a hypersurface an \emph{arrangement divisor}. Its closure in 
$\XX_\G^\bb$ is also a hypersurface, but need not be $\QQ$-Cartier. This is in a sense the reason that we have to consider modifications of the Baily-Borel compactification (and encounter them in algebraic geometry). 

We will write $\mathring{\XX}_\Hcal\subset \XX$ for the open arrangement complement $\XX-\cup_{H\in\Hcal}H$ and $\mathring{\LL}^\times_\Hcal$ for its preimage in $\LL^\times$. But when $\Hcal$ is understood, we may omit $\Hcal$ in the notation and simply write $\Delta$, $\mathring{\XX}$,  $\mathring{\LL}^\times,\dots$.
\\

\subsubsection*{Arrangement blowup.} We note here that there is a simple kind of blowup of $\widetilde\XX\to \XX$  in the analytic category with the property that the preimage of $\Delta$ is a normal crossing divisor: let $L(\Hcal)$ denote the collection of nonempty, proper subsets of $\XX$ obtained as an intersection of members of $\Hcal$. The minimal members of $L(\Hcal)$ are of course pairwise disjoint. Blow them up and repeat this construction for the strict transforms of the members of $\Hcal$. This process clearly stabilizes and the final result defines our  blowup $\widetilde\XX\to \XX$.
What is important for us  is no so much the fact that the preimage of $\Delta$ is a normal crossing divisor (note that this is even a nontrivial modification where $\Delta$ is already locally a genuine normal crossing), but that this modification admits an interesting blowdown. With this in mind we observe  that  $\widetilde\XX$ has a decomposition indexed by the members of $L(\Hcal)$ as follows.  Given $J\in L(\Hcal)$, then we have a natural trivialization of the normal bundle of $J$ in $\XX$ which is compatible with members of $\Hcal$ passing through $J$.  Denote the projectivized typical fiber by $E_J$ and denote by $\mathring{E}_J\subset E_J$ the corresponding arrangement complement. If  $\tilde J\to J$ is the blowup obtained by applying  the construction above to $\Hcal | J$, then the stratum associated with $J$ is naturally isomorphic to $\tilde J\times \mathring{E}_J$. The blowdown 
we alluded to is on this stratum given as projection on the second factor, so that the resulting analytic space has a decomposition $\mathring{\XX} \sqcup \bigsqcup_{J\in L(\Hcal)} \mathring{E}_J$.
We shall refer to $\widetilde\XX\to \XX$ as the \emph{arrangement blowup}. Its  blowdown 
$\widetilde\XX\to \XX$ (whose existence is proved in \cite{looij1} and \cite{looij2}) should be thought of as an alternate $\G$-equivariant extension, $\mathring{\XX}^{\alt}$  of $\mathring{\XX}$ in the analytic category; in particular, we get a modification $\XX_\G\leftarrow \widetilde\XX_\G\to \mathring{\XX}_\G^{\alt}$.

\section{Modular Examples}

\subsection*{Basics of GIT} Let $G$ be a complex reductive algebraic group and $H$ a finite dimensional complex  
representation of $G$. Then $G$ acts on the algebra $\CC[H]$ of regular functions on $H$. Here are two basic facts.
\begin{description}
\item[Finiteness] The graded $\CC$-algebra $\CC[H]^G$ is normal and finitely generated with generators in positive degree. So we have defined the weighted homogeneous (affine) cone $\spec (\CC[H]^G)$ whose base is the projective variety 
$\proj (\CC[H]^G)$. 
\item[Geometric interpretation] The closure of every $G$-orbit in $H$ contains a closed $G$-orbit and this $G$-orbit is unique. Any two distinct closed orbits are separated by a member of $\CC[H]^G$ and assigning to a closed $G$-orbit the corresponding point of $\spec \CC[H]^G$ sets up a bijection between the set of closed $G$-orbits on $H$ and the closed points of $\spec \CC[H]^G$ (we therefore write $G\bbs H$ for $\spec \CC[H]G$).
\end{description}
This also leads to an interpretation of $\proj \CC[H]^G$. The obvious map $H\to G\bbs H$ has as fiber through $0\in H$ the 
set of $v\in H$ with $0\in \overline{Gv}$. The complement,  called the \emph{semistable locus} and denoted $H^\sst$, then maps to the complement of the vertex of the cone $\spec \CC[H]^G$ and we get a bijection between the set of closed orbits in 
$\PP (H^{\sst})$ and the set of closed points of $\proj \CC[H]^G$. This is why we denote the latter sometimes by  $G\bbs\PP (H^\sst)$. Often these
separated quotients contain an genuine orbit space as an open-dense subset: the \emph{stable locus}  $H^\st\subset H^\sst$ consists of the $v\in H-\{ 0\}$ whose orbit is closed and for which the image of stabilizer $G_v$ in $\GL (H)$ is finite. Then $H^\st$ is Zariski open (but it need not be dense; it could even be empty while  $H^\sst$ is nonempty). Then the morphisms
$H^\st\to G\bbs H$ and  $\PP(H^\st)\to G\bbs\PP (H^\sst)$ are open with  image the $G$-orbit space of the source.
In case $H^\st$ is dense in $H^\sst$, we may regard $G\bbs\PP (H^\sst)$ as a projective compactification of the $G$-orbit space $G\bs H^\st$ of $H^\st$.
\\

Perhaps the most classical example is the following.  Given a complex vector space $U$ of dimension 2,  then for any integer $n>0$ we have the irreducible  $\SL (U)$-representation $H_n:=\sym^{n}(U^*)$. We regard $H_n$ as the space of homogeneous polynomials of degree $n$ on $U$ and $\PP(H_n)$ as the linear system of degree $n$ divisors on the projective line $\PP(U)$. Then $H_n^\sst$ resp.\ $H_n^\st$ parametrizes  the divisors which have no point of multiplicity $> n/2$ resp. $\ge n/2$. So when $n$ is odd, we have $H_n^\st=H_n^\sst$, so that $\SL (U)\bbs \PP(H_{n}^\sst)=\SL (U)\bs \PP(H_{n}^\st)$. When $n$ is even the difference $\PP(H_n^\sst)-\PP(H_n^\st)$ contains only one closed orbit, namely the one for which the divisor is given by two distinct points, each of multiplicity $n/2$ and so this will represent the unique point of $\SL (U)\bbs H_{n}^\sst -\SL (U)\bbs H_{n}^\st$.
In either case, the locus $H_n^\circ\subset H_n$ which defines reduced divisors ($n$-element subsets of  $\PP(U)$)  is the complement of the discriminant $D_n\subset H_n$, a hypersurface.

\begin{example}[Cyclic covers of the projective line following Deligne-Mostow]
Let us now focus on the case $n=12$ and write $H$ for $H_{12}$.
For $F\in H^\circ$ we have a cyclic cover of degree 6, $C_F\to \PP(U)$, ramified over the divisor $\Delta_{F}$ defined by $F$:
its homogeneous equation is $w^2=F(u)$, where it is understood that $u$ has degree 2 so that for every $\lambda\in \CC^\times$, $(\lambda^2u,\lambda w)$ and 
$(u,w)$ define the same point. This only depends on the image of  $F$ in $\PP(H)$, but this is not so for the differential $\omega_F$ we define next: fix a translation-invariant 2-form $\alpha$ on $U$.  Then  $\omega_F$ is a residue of the 2-form of $\alpha/w$ restricted to the surface defined by $w^6=F(u)$ at infinity. To be concrete: if $(u_0,u_1)$ is a coordinate pair for $U$, such that $\mu=du_0\wedge du_1$, then an affine equation for $C_F$ is $w'{}^6=F(1, u')$ with $u'=u_1/u_0$ and $w':=w/u_0^2$, and $\omega_F$ is in terms of these coordinates given by $du'/w'$.  The genus of $C_F$ is easily computed to be $25$.
The Galois group of the cover is the group $\mu_6$ of sixth roots of unity and acts on the $C_F$ via the $w$-coordinate according to the inverse of the (tautological) character $\chi: \mu_6\hookrightarrow \CC^\times$. So $\omega_F$ is an eigenvector for the character $\chi$.  The group $\mu_6$ has no invariants in $H^1(C_F,\CC)$ and the eigenspaces of $\mu_6$ with nontrivial characters all have dimension $10$.
Moreover, for $i=1,\dots , 5$, $\dim H^{1,0}(C_F,\CC)^{\chi^i}=2i-1$ with  $H^{1,0}(C_F,\CC)^{\chi}$ spanned by $\omega_F$. So $H^{1}(C_F,\CC)^\chi$ has signature $(1, 9)$. 

Now fix a unimodular symplectic lattice $V_\ZZ$ of genus 25 endowed with an action of $\mu_6$ such that  there exists a $\mu_6$-equivariant isomorphism $H^1(C_F,\ZZ)\cong V_\ZZ$ of symplectic lattices. Such an isomorphism will be unique up to an element of the centralizer of $\mu_6$ in  $\Sp (V_\ZZ)$. Notice that this centralizer (which we shall denote by $\Sp (V_\ZZ)_{\mu_6}$) contains a faithful image of $\mu_6$.  If $\G$ denotes the image of $\Sp (V_\ZZ)_{\mu_6}$ in $U(V_\CC^\chi)$, then $\omega_F$ defines an element of $(\LL_{V_\CC^\chi}^\times)_\G$ so that we have defined a map $H^\circ\to (\LL_{V_\CC^\chi}^\times)_\G$. This map is easily seen to be constant on the $\SL (U)$-orbits in $H^\circ$. Deligne and Mostow \cite{delmostow} show that this map extends to an isomorphism $\SL (U)\bbs H^\st\cong (\LL_{V_\CC^\chi}^\times)_\G$. One can show that this isomorphism extends from one of  $\SL (U)\bbs H^\sst$ onto $(\LL_{V_\CC^\chi}^\times)^\bb_\G$ so that the GIT compactification 
 $\SL (U)\bbs \PP(H^\sst)$ gets identified with the Baily-Borel compactification  $(\BB_{V_\CC^\chi})^\bb_\G$ (both are one-point compactifications). This also identifies their algebras of regular functions: the graded algebra of automorphic forms $A_\pt(\LL_{V_\CC^\chi}^\times)^\G$ gets identified with the algebra of invariants $\CC[H]^G$. 
\end{example}
 
\begin{remark}\label{rem:}
In particular, there is an automorphic form that defines the discriminant $D\subset H$. It would be interesting to see such an automorphic form written down. We expect this to have an infinite product expansion for the following reason: if two of the 12 points in $\PP(U)$ coalesce, then the curve $C_F$ is degenerate: it acquires an $A_5$-singularity with local equation $w^6=u^2$,  but the differential $\omega_F$ (locally like $du/w$) does not degenerate in the sense that it becomes a regular differential on the normalization of $C_F$ (this is clear if we write $u=\pm w^3$). This accounts for  the integral of $\omega_F$ over a vanishing cycle  to vanish. This implies that  near such a point $D$ maps to a hyperplane section of  $\LL_{V_\CC^\chi}^\times$.
\end{remark}

\begin{remark}\label{rem:}
We could also consider the other eigenspaces. This amounts  to passing to the intermediate Galois covers of the projective line: the $\mu_3$-cover and the hyperelliptic cover ramified in 12 points. The above isomorphism then leads to interesting morphisms from the $9$-dimensional ball quotient to a locally symmetric variety of type I$_{3,7}$ (of dimension 21) and to an arithmetic quotient of  a Siegel upper half space of genus 5 (of dimension 15).
\end{remark}

The theory of Deligne and Mostow provides many more examples, but the one discussed here is the one of highest dimension in their list. This list includes  the case of $\mu_4$-coverings of a projective line totally ramified in 8 points, which is related to the example that we discuss next.

\begin{example}[Quartic plane curves] Let $U$ be complex vector space of dimension 3 and let $H:=\sym^{4}(U^*)$. This is an irreducible representation of $\SL (U)$. So $\PP(H)$ is the linear system of degree 4 divisors on the projective plane $\PP(U)$. We denote the divisor associated with $F\in H-\{ 0\}$ by $C_F$ and we let $H^\circ\subset H$ be the set of $F\in H$ for which $C_F$ is a smooth quartic curve.  This is the complement of the discriminant hypersurface $D\subset H$.
It is classical fact that $F\in H^\st$ if and only if $C_F$ is a reduced curve whose singularities are only nodes or cusps and that any closed orbit of $H^\sst-H^\st$ is representable by $(u_1u_2-u_0^2)(su_1u_2-tu_0^2)$ for some $(s,t)\in\CC^2$ with $s\not=0$
(hence defines the sum of  two conics, one of which is smooth, which meet at two points with multiplicity $\ge 2$). This includes  the case of a smooth conic with multiplicity 2. The orbit space $\SL(U)\bs \PP(H^\circ)$ has a simple modular interpretation: it is the moduli space of nonhyperelliptic curves of genus 3 (for each such curve is canonically embedded  in a projective plane as a quartic curve). 

With every $F\in H^\circ$ we associate the $\mu_4$-cover  $S_F\to\PP (U)$ which totally ramifies along $C_F$, i.e., the  smooth quartic surface in $\PP(U\oplus \CC)$ defined by $w^4=F(u)$; indeed, this is a polarized K3 surface of degree 4 with 
$\mu_4$-action. The same argument as in the previous example shows that if $\alpha$ is a translation $3$-form on $U$, then the residue of  $\alpha/w^3$ defines a nowhere zero 2-form $\omega_F$ on $S_F$. We let $\mu_4$ act  on $w$ with the tautological character 
$\chi :\mu_4\hookrightarrow \CC^\times$ so that $\omega_F\in H^{2,0}(S_F)^\chi$. In fact, $\omega_F$ generates $H^{2,0}(S_F)$. 
The abelian group $H^2(S_F,\ZZ)$ is free abelian of rank $22$ and the intersection pairing  $H^2(S_F,\ZZ)\times H^2(S_F,\ZZ)\to \ZZ$ is even unimodular of signature $(3,19)$.
The $\mu_4$-fixed point subgroup in $H^2(S_F,\ZZ)$ is the image of $H^2(\PP(U),\ZZ)\to H^2(S_F)$ and contains the polarization of $S_F$. The other character spaces in $H^2(S_F,\CC)$ all have dimension 7.
Since  $H^{2,0}(S_F,\CC)^\chi=H^{2,0}(S_F,\CC)$ is spanned  by $\omega_F$, it follows that the orthogonal complement of $H^{2,0}(S_F,\CC)^\chi$ (relative to the Hermitian form) in  $H^2(S_F,\CC)^\chi$ is primitive cohomology of type $(1,1)$ and of dimension $6$. In particular,  $H^{2,0}(S_F,\CC)^\chi$ has signature $(1,6)$.

Fix a unimodular even lattice $\Lambda$ of signature $(3,19)$ endowed with an action of $\mu_4$ such that  there exists a $\mu_4$-equivariant isomorphism $\phi: H^2(S_F,\ZZ)\cong \Lambda$ of lattices. Let us write $(W,h)$ for the Hermitian vector space $\Lambda_\CC^\chi$ and let $\LL^\times$ be the set of $w\in W$ with $h(w,w)>0$. Then $\phi$  will be unique up to an element of $\Orth(\Lambda)_{\mu_4}$. If $\G$ denotes the image of this group in $U(W)$, then $\omega_F$ defines an element of $\LL^\times_\G$ so that we have defined a map 
$H^\circ\to \LL^\times_\G$. This map is  constant on the $\SL (U)$-orbits in $H^\circ$. Kond\=o has shown \cite{kondo} that this map extends to an open embedding $\SL (U)\bbs H^\st\hookrightarrow\LL^\times_\G$. 

The situation is however not as nice as in the previous example: the map is not surjective and, related to this, does not extend to $\SL (U)\bbs H^\sst$. In fact, its image is the complement of an arrangement divisor. This can be explained by the fact that we miss out some K3 surfaces of degree 4 with $\mu_4$-action. It is also related to the fact that we miss out some of the genus 3 curves, namely the hyperelliptic ones. This locus is represented in $\LL^\times_\G$ as an arrangement divisor $\Delta_\G$ so that we have an 
isomorphism  $\SL (U)\bbs H^\st\cong (\mathring{\LL}^\times)_\G$. The divisor $\Delta_\G$ is irreducible: 
the group $\G$ acts transitively on $\Hcal$,  and its normalization, the quotient of a member $H\in\Hcal$ by its $\G$-stabilizer, 
is a copy the Deligne-Mostow ball quotient for the pair $\mu_4$-covers of projective line totally ramified in 8 points.  This locus is not visible in $\SL (W)\bbs \PP(H^\sst)$, for  $\PP(H^\sst)-\PP(H^\st)$  is just a singleton (represented by the conic with multiplicity 2). Since this difference is codimension one, $\CC[H]^{\SL (U)}$ is also the algebra of regular functions on
$\SL (U)\bbs H^\st$. So via the isomorphism  $\SL (U)\bbs H^\st\cong (\mathring{\LL}^\times)_\G$ the graded algebra $\CC[H]^{\SL (U)}$ is reproduced as the graded algebra of $\G$-invariant analytic functions on $ (\mathring{\LL}^\times)_\G$ that are sums of homogeneous functions. There is no growth condition here and  these functions are automatically meromorphic on  $\LL^\times$. From the  ball quotient perspective, it is quite a surprise that this algebra is finitely generated and has positive degree generators!
\end{example}

\begin{example}[Double covers of sextic curves]
Let $U$ be complex vector space of dimension 3 and let $H:=\sym^{6}(U^*)$ be a $\SL (U)$-representation. So $\PP(H)$ is the linear system of degree 6 divisors on the projective plane $\PP(W)$. As above, we  denote by $H^\circ\subset H$ be the set of $F\in H$ for which $C_F$ is a smooth sextic curve. According to Jayant Shah \cite{shah}, $H^\st$ contains $H^\circ$ and allows $C_F$ to have simple singularities in the sense of Arnol'd. The closed orbits in $\PP(H^\sst-H^\st)$ make up a longer list and come in families, but two of them deserve special mention: the closed orbit represented by a smooth conic of multiplicity 3 and the 
the closed orbit represented by a coordinate triangle of multiplicity 2. We proceed as before: we associate to $F\in H^\circ$ the double cover $S_F\to \PP(W)$ which totally ramifies along $C_F$: it has the weighted homogeneous equation $w^2=F(u)$ (where $\deg w=3$) and $\omega_F:=\alpha/w$ defines a nowhere zero regular  2-form on $S_F$. The surface $S_F$ is a K3 surface. The obvious  involution of $S_F$ acts on $H^2(S_F,\ZZ)$  with the image of  $H^2(\PP(U),\ZZ)$  in $H^2(S_F,\ZZ)$  (a copy of $\ZZ$) as its fixed point set. The sublattice  $H^2(S_F,\ZZ)^-\subset H^2(S_F,\ZZ)$ on which this involution acts as minus the identity is nondegenerate and of signature $(2,19)$. Notice that its complexification contains $\omega_F$.

We therefore  take the same lattice $(\Lambda, s)$ of signature $(3,19)$ as above, but now endowed with an involution $\iota$ such that  there exists an isomorphism $\phi:H^2(S_F,\ZZ)\cong \Lambda$ of lattices with involution. We  denote by
$V\subset \Lambda\otimes\RR$  the subspace on which $\iota$ acts as minus the identity (it has signature $(2,19$) and  write  $\LL^\times$ for $\LL^\times_V$ (a connected component of the set $v\in V_\CC$ for which $s_\CC(v, v)=0$ and $s_\CC(v,\overline v)>0$).
We ask that the isomorphism $\phi$  maps $\omega_F$ to $V^+$ and then $\phi$ will be unique up to an element of $\G:=\Orth^+(\Lambda)_{\iota}\subset \Orth (V_\ZZ)$. The above construction produces a map 
$\SL(U)\bs H^\circ\to \LL^\times_\G$. This map is  constant on the $\SL (U)$-orbits in $H^\circ$ and J.~Shah has shown that it extends to an open embedding $\SL (U)\bbs H^\st\hookrightarrow \LL^\times_\G$. But he also observed that as in the previous case its map is not onto: the image is the complement of an irreducible arrangement divisor $\Delta_\Hcal$ and the map does not extend to $\SL (U)\bbs H^\sst$. The explanation is similar: this divisor parametrizes the K3 surfaces of degree 2 that we missed, namely the hyperelliptic ones. Since $\SL (U)\bbs H^\sst-\SL (U)\bbs H^\st$ is of dimension 2 (hence of codimension $>1$ in $\SL (U)\bbs H^\sst$, $\CC[H]^{\SL (U)}$ can be understood as the graded algebra of $\G$-invariant analytic functions on $ \mathring{\LL}^\times$ that are sums of homogeneous functions (these functions are automatically meromorphic on  $\LL^\times$). Again there seems reason to be pleasantly surprised.
\end{example}

\subsection*{A potential example: Allcock's 13-ball} Consider an even unimodular lattice $V_\ZZ$ of signature $(2, 26)$ endowed with an action of $\mu_3\subset\Orth (V_\ZZ)$ which leaves no nonzero vector fixed. The pair $(V_\ZZ, \mu_3\subset\Orth (V_\ZZ))$ is unique up to isomorphism. If $\chi :\mu_3\subset\CC^\times$ denotes the tautological character, then $V^\chi_\CC$ has signature $(1,13)$ and so we have a $13$-dimensional ball 
$\BB_{V^\chi_\CC}$. Allcock  \cite{allcock} makes a number of intriguing conjectures about the orbit space of this ball relative to  the group $\Orth (V_\ZZ)_{\mu_3}$ and  suspects that it has a modular interpretation (that is why we call it a potential example).


\section{A Baily-Borel package for arithmetic arrangement complements}
In some of the preceding two examples we found that an algebra of invariants can be understood as an algebra $\G$-invariant functions on an arrangement complement. These  $\G$-invariant functions are meromorphic, when considered as functions on $\LL^\times$: they are meromorphic automorphic forms. This cannot be just a coincidence: it suggests that there should be an Baily-Borel package in that setting. This is indeed the case. Let us first focus on the case of a complex ball. We shall always assume  that $\G$ is \emph{neat} in the sense of Borel: this means that  for every $\g\in\G$, the subgroup of $\CC^\times$ generated by the eigenvalues of $\G$ is
torsion free.  This is fairly inessential assumption, as this  can always be satisfied by passing to a subgroup of finite index. We are going to construct a chain of normal proper varieties of the following type
\[
\XX_\G^\bb\leftarrow  \XX_\G^\Hcal\leftarrow\widetilde{\XX}_\G^\Hcal \to \mathring{\XX}_\G^\bb,
\]
where $\XX_\G^\bb$ is the classical Baily-Borel compactification, the first morphism is a modification of that compactification over its boundary (and so an isomorphism over $\XX_\G$), which then makes it possible to
carry out  the arrangement blowup (this is subsequent morphism) and  a blowdown across the boundary. So  the  the last morphism is a blowdown that should be thought of as extending the formation of the alternate extension $\mathring{\XX}^\alt_\G$ of $\mathring{\XX}_\G$ across boundary, as it does not affect $\mathring{\XX}_\G$.  But as the notation suggests, we like to think of the contracted space  as a Baily-Borel compactification of  that variety. We will refer to this as the \emph{modification chain} associated with $\Hcal$. 

\subsection*{Arithmetic arrangements on complex balls} Let $(W,h)$ be a hermitian vector space of Lorentz signature $(1,n)$ with $n\ge 1$  and  $\G\in U(W)$ the image of an arithmetic group acting on $W$ (so this presupposes that we have a CM-subfield $K\subset \CC$ and a $K$-form $W_K\subset W$ such that $h|W_K\times W_k$ takes its values in $K$, and $\G\subset \U(W_K)$).  Let be given a $\G$-arrangement  $\Hcal$ on $W$: every $H\in \Hcal$ is a hyperplane of $W$  defined over $K$ and of  signature $(1, n-1)$, and the collection $\Hcal$ is a finite union of $\G$-orbits. We denote by $\mathring{\BB}$ the associated arrangement complement and 
 $\mathring{\LL}^\times$ its preimage in ${\LL}^\times$.
 
Denote by $\Jcal_\Hcal$ the  collection of subspaces $J\subset W$ that are not positive  for $h$ and can be written as an intersection of members of the collection of hyperplanes $\Hcal\cup\Ical^\perp$. Here $\Ical^\perp$ is the set of $I^\perp$ with $I\subset W$ an isotropic line. Notice that $W$ belongs to $\Jcal_\Hcal$, as that corresponds to the empty intersection, but that $\{0\}$ does not, because $\{ 0\}$ is positive definite. We will describe  for every member the image of $\mathring{\LL}^\times$  under the projection $\pi_J: W\to W/J$.

Suppose first that $J\in \Jcal_\Hcal$ is nondegenerate: so $J$ Lorentzian and $J^\perp$ is negative definite. Then $J^\perp$ maps isomorphically to $W/J$ and the $\G$-stabilizer of $J$ will act on $W/J$ through a finite group. The image of  $\LL^\times$ in   $W/J$ is all of  $W/J$, whereas the image of $\mathring{\LL}^\times$ is the complement in $W/J$ of the union of the hyperplanes $H/J$  with $H\in\Hcal$ and $H\supset J$ (there are only finitely many such).  

Next consider the case when $J\in \Jcal_\Hcal$ is degenerate. Then $I:=J\cap J^\perp$ is an isotropic line and we have $I\subset J\subset I^\perp$ and  $I\subset J^\perp\subset I^\perp$. We have a factorization
\[
\begin{CD}
\LL^\times\subset W@>{\pi_I}>>W/I @>{\pi^I_J}>> W/J @>{\pi^J_{I^\perp}}>> W/I^\perp.
\end{CD}
\]
Since $I^\perp$ is negative semidefinite, we have $\LL^\times\cap I^\perp=\emptyset$. In fact, the image 
of $\LL^\times$ in  $W/I$ is  complement of the hyperplane $I^\perp/I$, its 
image  in $W/J$ is the complement of the hyperplane $I^\perp/J$ and its image in the one-dimensional $W/I^\perp$ the complement of $0$.  So these images all have the structure of $\CC^\times$-bundles over affine spaces. On the other hand, the first projection $\pi_I :\LL^\times\to W/I$ is a bundle of upper half planes.

The $\G$-stabilizer $\G_I$ of $I$ is a Heisenberg group. Its center (which is infinite cyclic) acts trivially on $W/I$, and as a discrete  translation group in the fibers of $\pi_I :\LL^\times\to W/I$ and the (abelian) quotient by this center acts faithfully on every fiber  of $W/I\to W/I^\perp$ as a full lattice in  the vector group $\Hom(W/I^\perp, I^\perp/I)$. The $\G_I$-orbit space of $W/I$  is a torsor  over $(W/I^\perp) -\{ 0\}$ (which is simply a copy of $\CC^\times$) with respect to an abelian variety that is isogenous to a product of elliptic curves. The  $\G_I$-orbit space of $\LL^\times$  has  the structure of a punctured disk bundle over this torsor. We may fill in the zero section of that punctured disk bundle as to get a disk bundle and then we have produced a local extension of the orbit space $\G\bs \LL^\times$ of Mumford type. This zero section can be contracted along the projection onto $W/I^\perp -\{ 0\}$, and then we have a local extension of Baily-Borel type. But we can also interpolate between the two and contract along the projection onto $W/J$ (before or after passing to $\G$-orbit spaces).

The $\G$-stabilizer of $J$ is contained in the $\G$-stabilizer of $I$ and it intersects the said Heisenberg group in an Heisenberg subgroup (with the same center), with the abelian quotient now acting  on any fiber of $W/J\to W/I^\perp$ as a full lattice of translations in $\Hom(W/I^\perp, J/I)$. So the $\G_J$-orbit space of $\pi_J(\LL^\times)$ is naturally a finite quotient of a $\CC^\times$-bundle an abelian  torsor. If we turn our attention to  
$\pi_J(\mathring{\LL}^\times)$, then a finite number of abelian subtorsors of codimension one must be left out: these 
are defined by  the $H\in \Hcal$ which contain $J$. Any such $H$ defines a hyperplane $H/J$ in $W/J$ and the 
stabilizer $\G_J$ has only finitely many orbits in this collection such hyperplanes.

The disjoint union
\[
(\mathring{\LL}^\times)^\bb:= \mathring{\LL}^\times\sqcup \bigsqcup_{J\in\Jcal_\Hcal} \pi_J(\mathring{\LL}^\times),
\] 
can be endowed with a Satake type of topology that is invariant under $\G\times\CC^\times$. We then find that if $\G$ is sufficiently small (in the sense that no eigenvalue $\not=1$ of an element of $\G$ is of finite order), then the orbit space $(\mathring{\LL}^\times)^\bb_\G$ is an extension of  $\mathring{\LL}^\times$ by a finite number strata, each of which is either  a linear arrangement complement or an affine arrangement complement. With $J=V$ is associated a singleton stratum. If we remove that singleton and divide by $\CC^\times$, we get a space denoted $\mathring{\BB}^\bb$.

If we are lucky and all these strata have codimension $>1$ (which means that no member of $\Jcal$ has dimension 1, then a Koecher principle applies and we have an easily stated  Baily-Borel package:
\begin{description}
\item[Finiteness] If $A_d (\mathring{\LL}^\times)$ denotes the space of holomorphic functions on $\mathring{\LL}^\times$ that are homogeneous of degree $-d$, then the graded $\CC$-algebra $A_\pt (\mathring{\LL}^\times)^\G$ is finitely generated with generators in positive degree. So we have defined the weighted homogeneous (affine) cone 
$\spec A_\pt (\mathring{\LL}^\times)^\G$ whose base is the projective variety $\proj (A_\pt (\mathring{\LL}^\times)^\G)$. 
\item[Separation and Topology] The underlying topological spaces (for the Hausdorff topology) are naturally  identified with $(\mathring{\LL}^\times)^\bb_\G$ resp.\ $\mathring{\BB}^\bb_\G$ and via these identifications,
$\spec A_\pt (\mathring{\LL}^\times)^\G$ and its projectivization $\proj (A_\pt (\mathring{\LL}^\times)^\G)$ acquire a partition into locally closed subvarieties (the former invariant under $\CC^\times$).
\end{description}
When there are strata of codimension one, this still holds if we  impose certain growth conditions. We think of $A_\pt (\mathring{\LL}^\times)^\G$ as an algebra of meromorphic $\G$-automorphic forms (we allow poles along the hyperplane sections indexed by $\Hcal$).

For $\Hcal=\emptyset$, we get the usual Baily-Borel compactification $\BB^\bb_\G$ of $\BB_\G$. There is a  natural compactification of $\BB_\G$ which dominates $\BB^\bb_\G$ and is the first step up towards a description of  $\mathring{\BB}^\bb_\G$: we then let  for every isotropic line  $I$,  $J_I$ be the intersection of $I^\perp$ and the $H\in\Hcal$ with $H\supset I$. We may then form the disjoint union
\[
(\LL^\times)^\Hcal:= \LL^\times \sqcup \bigsqcup_{I\in\Ical} \pi_{J_I}({\LL}^\times),
\] 
and its $\CC^\times$-orbit space $\BB^\Hcal:= \PP((\LL^\times)^\Hcal)$. There is an obvious $\G$-equivariant map from $\BB^\Hcal$ to $\BB^\bb$. The arrangement blowup and blowdown  as described in Section \ref{section:bbpackage} can now be carried out on $\BB^\Hcal$ without much change, but the setting is that within a category of a locally ringed spaces. This yields the modification chain
\[
\BB^\bb_\G \leftarrow \BB^\Hcal_\G \leftarrow\widetilde{\BB}^\Hcal_\G \to\mathring{\BB}^\bb_\G. 
\]
If we recall that $\BB^\bb_\G $ is obtained from $\BB_\G $ by adding a finite set of cusps, then the first morphism is a modification of $\BB^\bb_\G$ over the cusps. The exceptional locus of this modification need not be codimension one. But the  closure of the arrangement $\Delta_\G$  in  $\BB^\bb_\G$ has the pleasant property that becomes a Cartier divisor whose strata we can blow up and  down in a specific manner as to form $\mathring{\BB}^\bb_\G$.

\subsection*{A revisit of the moduli space of quartic curves}
Let us see how this works out in the example of a quartic curve.
Recall that Kond\=o's  theorem states that $\SL(U)\bs H^\st$ maps isomorphically onto an arrangement complement $(\LL^\times_{\Hcal})_\G$. This assertion can now be amplified: since isomorphism 
gives rise to an isomorphism of $\CC$-algebras $\CC[U]^H\cong A_\pt (\mathring{\LL}^\times)^\G$, this 
isomorphism must extend to an isomorphism of $\SL(U)\bbs H^\sst$ onto $(\LL^\times_{\Hcal})^\bb_\G$.
But the small blowup of the Baily-Borel compactification also appears here:  we have a closed strictly semi-stable orbit  in $\PP(H^\sst)$ of conics with multiplicity 2. This orbit hides away from us the genus 3 curves that cannot be canonically embedded, namely  the hyperelliptic curves. But they are unlocked by first blowing up this orbit in
$\PP(H^\sst)$ to get $\widetilde{\PP}(H^\sst)$ and then form the universal separated $\SL (U)$-orbit space (this replaces the point representing the orbit in question by the coarse moduli space of pairs consisting of a nonsingular conic with an effective degree 8 divisor that is Hilbert-Mumford stable). The resulting  space projects to the Baily-Borel compactification 
$\BB^\bb_\G$. The latter has just one cusp and the fiber over that cusp is a rational curve, parametrized by the family of quartics $[s:t]\in\PP^1\mapsto [(u_1u_2-u_0^2)(su_1u_2-tu_0^2)]$. For $[s:t]=[1:1]$ we land in the point represented by the double conic $[(u_1u_2-u_0^2)^2]$ endowed with the (strictly semistable) divisor $4(p_1)+4(p_2)$, where $p_1=[0:1:0]$ and $p_2=[0:0:1]$. We may understand the projection
\[
 \SL (U)\bbs\widetilde{\PP}(H^\sst)\to \SL (U)\bbs\PP(H^\sst)
\]
as the morphism $\BB^\Hcal_\G \to \mathring{\BB}^\bb_\G$. The latter is a morphism indeed: the closure of the arrangement divisor in $\BB^\Hcal_\G$ needs no blow-up (we have $\widetilde\BB^\Hcal_\G=\BB^\Hcal_\G$) in order to contract it. This is because distinct members of $\Hcal$ do  not meet in $\BB$.

\subsection*{Arithmetic arrangements for type IV domains}
Let $(V,s)$ be a nondegenerate symmetric bilinear form of signature $(2,n)$ defined over $\QQ$ with $n\ge 1$, $\G\subset \Orth^+(V)$ an arithmetic group and $\Hcal$ a $\G$-arrangement in $V$. We write $\LL^\times$ for $\LL^\times_V$, a connected component of the set of  $v\in V_\CC$ with $s_\CC(v,v)=0$ and $s_\CC(v,\overline v)>0$ and $\DD\subset \PP(V)$ for its projectivization. 

We are going to define an indexed collection $\{J_\alpha\}_{\alpha\in \Jcal_\Hcal}$ of nonpositive subspaces $J\subset V$ defined over $\QQ$. The indexing is an essential part of the data that cannot be ignored, because the map $\alpha\in \Jcal_\Hcal\mapsto J_\alpha$ will in general not be injective. (The reason is that the strata are often defined by torus embeddings and are therefore indexed by locally polyhedral cones rather than by linear subspaces.) But apart from this,  once this collection is defined, we can define $(\mathring{\LL}^\times)^\bb$ as in the ball quotient case and assert that the Baily-Borel package holds verbatim.

The subspaces $J_\alpha$ in question will all have the property that $J_\alpha^\perp$  is negative semidefinite. We break them up into  the three classes defined by the dimension of their radical 
$J_\alpha\cap J^\perp_\alpha$, which can be $0$, $1$ or $2$:  $\Jcal_\Hcal=\sqcup_{i=0}^2 \Jcal_\Hcal^{(i)}$.  We will describe these collections in detail below. But even without this precise information, we can already say what the 
modification chain is like. In  the special case when $\Hcal$ is empty,  $\Jcal=\Jcal_\Hcal$ is simply the collection of $I^\perp$, where $I$ is an isotropic subspace defined over $\QQ$. We also need a subcollection 
$\Jcal^\infty_\Hcal\subset \Jcal_\Hcal$: this will be the union of $\Jcal_\Hcal^{(1)}$ and the minimal members of 
$\Jcal_\Hcal^{(2)}$. Then we can form
\begin{align*}
\DD^\bb: &= \DD \sqcup\bigsqcup_{J\in \Jcal} \PP(\pi_J)(\DD),\\
\DD^\Hcal: &= \DD \sqcup\bigsqcup_{J\in \Jcal^\infty_\Hcal} \PP(\pi_J)(\DD),\\
\DD_\Hcal^\bb: &= \DD_\Hcal \sqcup\bigsqcup_{J\in \Jcal_\Hcal} \PP(\pi_J)(\DD_\Hcal).\\
\end{align*}
At this point, these are just disjoint unions of complex manifolds on which $\G$ acts. We have a projection
\[
\DD^\Hcal\to \DD^\bb
\]
defined as follows.  For $J\in \Jcal_\Hcal^\infty$, we note that $J+J^\perp$ is in $\Jcal$  and on the associated stratum  this map is simply be the projection $\PP(\pi_J)(\DD)\to \PP(\pi_{J+J^\perp})(\DD)$. 

So let us then define $\Jcal_\Hcal$. Our discussion will also clarify what the topology on these  spaces must be. Since the collections $\Jcal_\Hcal^{(0)}$ and $\Jcal_\Hcal^{(2)}$ are easiest to treat (they are self-indexing) we describe them first.
\\

\subsubsection*{The collection $\Jcal_\Hcal^{(0)}$.}
The collection  $\Jcal_\Hcal^{(0)}$ consists of all the linear subspaces $J\subset V$   for which $J^\perp$ is negative definite and  that appear as an intersection of members of $\Hcal$ (including $V$). We regard this as a self-indexing collection. For any such $J$, $V/J $  can be identified with $J^\perp$. The $H\in\Hcal$ containing $J$ define a finite arrangement in $V/J$ and the $\G$-stabilizer $\G_J$ of $J$ acts on $V/J$ through a finite group.  The projection $\pi_J: V_\CC\to V_\CC/J_\CC$ is surjective when restricted to $\LL^\times$ and the image of  $\mathring{\LL}^\times$ is the complement of the complexified arrangement in $V_\CC/J_\CC$. The linear section
$\DD\cap \PP(J)$  is a totally geodesic complex submanifold of type IV on which $\G_J$ properly as an arithmetic group.
\\

\subsubsection*{The collection $\Jcal_\Hcal^{(2)}$.} 
Similarly,  the collection  $\Jcal_\Hcal^{(2)}$ consists of all the intersections  $J\subset V$ of members of $\Hcal$ with radical $I:=J^\perp\cap J$ is of dimension 2; this too is a self-indexing collection. The situation here is  not unlike that of a cusp in the ball quotient case. We have $I\subset J\subset I^\perp$ and  a factorization
\[
\begin{CD}
\LL^\times\subset V_\CC@>{\pi_I}>>V_\CC/I_\CC @>{\pi^I_J}>> V_\CC/J_\CC@>{\pi^J_{I^\perp}}>> V_\CC/I^\perp_\CC\cong \Hom( I, \CC).
\end{CD}
\]
We first concentrate on  the structure that is imposed by $I$ and then return to $J$.
As we observed before, the plane $I$ is naturally oriented by our choice of the connected component $\LL^\times$ of $V^+$ and the subset $\pi_{I^\perp}(\LL^\times)$ is the full set $\iso^+(I, \CC)$ of orientation preserving  real  isomorphisms $\phi: I\cong \CC$.  Its projectivization 
$\PP (\iso^+(I, \CC))$ is  a  half plane in the Riemann sphere  $\PP(\Hom( I, \CC))$ and parametrizes all the complex  structures on $I$ compatible with its orientation. The arrangement complement $\mathring{\LL}^\times$ maps still onto  $\iso^+(I, \CC)$. 
The $I$-stabilizer of $\G$,  $\G_I$,  acts on this image through an arithmetic subgroup $\G (I)$ of $\SL (I)$ and so the orbit space of $\pi_{I^\perp}(\LL^\times)$ is a $\CC^\times$-bundle over modular curve; we denote that orbit space by $\iso^+(I, \CC)_{\G(I)}$. 

Let us now concentrate on the action of the subgroup $Z_\G(I)\subset \G_I$ which fixes $I$ pointwise. This group preserves every fiber of  $\pi_{I^\perp}: \LL^\times\to \iso^+(I, \CC)$. We have a kind of dual situation regarding the projection $\LL^\times\to \pi_I(\LL^\times)$: this  is a bundle of upper half planes. The intermediate projection $\pi_I(\LL^\times)\to\iso^+(I, \CC)$ has the structure  of a real affine bundle whose structural group the vector group $\Hom(I^\perp/I,I)$.  But if we specify a  $\phi\in \iso^+(I, \CC)$, then this vector group acquires a complex structure via its identification with  $\Hom(I^\perp/I,\CC)$ and hence the same is true for the fiber $\pi_I(\LL^\times_\phi)$ of $\pi_I(\LL^\times)\to\iso^+(I, \CC)$ over $\phi$. The group $Z_\G(I)$ is a Heisenberg group that acts faithfully on $\LL^\times_\phi$ as follows. Its (infinite cyclic) center can be identified with
$Z_\G(I^\perp)$: it acts trivially
on the affine space $\pi_I(\LL^\times_\phi)$ (observe that $V/I$ may be identified with the dual of $I^\perp$) and as a translation group in each upper half plane fiber of
$\LL^\times_\phi\to\pi_I(\LL^\times_\phi)$. The quotient $Z_\G (I)/ Z_\G (I^\perp)$ by the center  can be identified with a lattice in  $\Hom(I^\perp/I,I)$ and  acts via this lattice on  $\pi_I(\LL^\times_\phi)$. 

So the  $Z_\G (I)$-orbit space of $\LL^\times$ has the structure of a bundle over $\iso^+(I, \CC)$ whose typical  fiber is a punctured disk bundle over an abelian variety isogenous to a product of elliptic curves.  Filling in the punctures defines a divisorial extension of this orbit space as it appears in Mumford's toroidal compactification theory as exposed in \cite{amrt}. This boundary divisor is a copy of the $Z_\G (I)$-orbit space of $\pi_I(\LL^\times)$ and  has therefore the structure of an abelian torsor over  $\iso^+(I, \CC)$. The resulting disk  bundle sits in an line bundle 
that is anti-ample relative to the projection onto $\iso^+(I, \CC)$. This implies that it can in fact be contracted analytically along the projection onto  $\iso^+(I, \CC)$. This replaces the toric divisor  by a copy of $\iso^+(I, \CC)$ (and so by a copy of a $\CC^\times$-bundle over a modular curve if we pass to the $\G_I$-orbit space):  we get a `minimal extension' of Baily-Borel type:   

The intermediate space $I\subset J\subset I^\perp$  gives rise to an intermediate extension: the projection 
$\pi_{I}(\LL^\times)\to \pi_J(\LL^\times)$ becomes after passage to $\G_I$-orbit spaces a projection 
of abelian torsors over  $\G (I)\bs\iso^+(I, \CC)$. The anti-ampleness property mentioned above enables us to accomplish this  as an analytic contraction inside the Mumford extension.  The resulting analytic  extension of the $\G_I$-orbit space of $\LL^\times$ then adds a copy of the  $\G_J$-orbit space of $ \pi_J(\LL^\times)$. The projection
$\pi_J(\LL^\times_{\Hcal})$ is the complement of the union of finitely many abelian subtorsors in $\pi_J(\LL^\times)$
of codimension one.  Indeed, the union $\G_I$-orbit space of $\LL^\times_{\Hcal}$ the  $\G_I$-orbit space of $ \pi_J(\LL^\times)$ defines an open subset of the extension just described.

The biggest extension obtained in this manner is when we take  $J$ to be  the intersection $J_I$ of $I^\perp$  with all the members of $\Hcal$ which contain $I$. If the latter set is empty, then $J_I=I^\perp$ and we obtain an extension of Baily-Borel type. If $J_I=I$, we get the Mumford toroidal extension. The extension defined by $J_I$ has the property that the closure for any $H\in \Hcal$, the closure of  $\LL^\times|H$ in this extension is a $\QQ$-Cartier divisor (this closure meets the added stratum only if $H\supset I$).
\\

\subsubsection*{The collection $\Jcal_\Hcal^{(1)}$.} 
The remaining elements $\Jcal_\Hcal$ parametrize  intersections  of members of $\Hcal$ with one-dimensional radical, but here the indexing is more subtle, because it involves the data of a `semi-toric embedding'.
Let us first fix  an isotropic subspace $I\subset V$ of dimension 1 that is defined over $\QQ$. So the form $s$ induces on $I^\perp/I$ a nondegenerate symmetric bilinear $s_I$ form of signature $(1,n-1)$ such that the pair $(I^\perp/I, s_I)$ is also defined over $\QQ$. The set of 
$v\in I^\perp/I$ with $s_I(v,v)>0$ has two connected components, each being a quadratic (Lobatchevski) cone, but there is a priori no preferred one. Similarly, $I$ has no preferred orientation. But the choice of the connected component $V^+=\LL^\times$ makes that  these choices are no longer independent, for  we will see that we have in $I\otimes  I^\perp/I$ a naturally defined quadratic cone $C$. 

Given $\omega\in \LL^\times$, then its image $\pi_{I^\perp}(\omega)$  in the complex line $(V/I^\perp)_\CC\cong \Hom_\RR(I, \CC)$ is nonzero. Hence so its image $\pi_I(\omega)$  in $(V/I)_\CC$. It thus determines a linear section of the projection  $(V/I)_\CC\to (V/I^\perp)_\CC$ which takes $\pi_{I^\perp}(\omega)$  to $\pi_I(\omega)$. This section only depends on the $\CC^\times$-orbit of $\omega$ so that we have in fact defined map from $\DD=\PP(\LL^\times)$ to the complexification $A(\CC)$ of the affine space $A$ of sections of $V/I \to V/I^\perp$. This  affine space  has as its translation space $\Hom (V/I^\perp, I^\perp/I)\cong I\otimes I^\perp/I$ so that the  complexification $A(\CC)$ of $A$ is naturally identified with $A+ \sqrt{-1}I\otimes I^\perp/I$. Let us just write $N$ for $I\otimes I^\perp/I$. The map just described is an open embedding of  $\DD$ in $A_I(\CC)$ and gives a realization of $\DD$ as a \emph{tube domain}:  
$s_I$ defines on  $N$ a naturally defined symmetric bilinear form 
$\tilde s_I:N \times N\to I\otimes I$. The target is a one dimensional vector space that has a distinguished positive half-line and the set of $u\in N$ with $\tilde s_I(u,u)$ in this half line 
consists of two opposite quadratic cones. One of these cones, denote it by $C$, has the property that $\DD$ maps isomorphically onto the open subset $A+\sqrt{-1}C$. This  open subset $A+ \sqrt{-1}C$ of
$A(\CC)$ is a classical tube domain. Our assumption regarding $\G$ implies that $\G_I$ leaves $I$ pointwise fixed. The group $\G_I$ is now realized as a group of motions in $A(\CC)$ which preserves this tube domain. In fact, the subgroup of $\G_I$ that acts as the identity on $I^\perp/I$ is abelian and is realized in $N$ as a lattice and acts as such (as a translation group) on $A$. We therefore denote that group by $N_\ZZ$.

The quotient $\G_I/N_\ZZ $ can be regarded as a subgroup of $\GL(N)$ or of $\GL(I^\perp/I)$. This subgroup (that we denote $\G(N)$) respects $C$ and the hyperbolic form attached to it; it is in fact arithmetic in that group. The $N_\ZZ$-orbit space of  $A(\CC)$, $T= A(\CC)_{N_\ZZ}$, is torsor for the algebraic torus $\TT:=N_\ZZ\otimes\CC^\times$. Tensoring the homomorphism $z\in \CC^\times\mapsto \log |z|\in \RR$  with $N_\ZZ$ gives a homomorphism $\Im: \TT\to N$ with kernel  the compact torus $N_\ZZ\otimes U(1)$ and we then can identify $\Tcal:=\Im^{-1}(C)$ with $\DD_{N_\ZZ}$.   We are going to apply the theory of semi-toric embeddings to this situation. 

Let $\Hcal_I$ denote the collection of $H\in\Hcal$ containing $I$. This a finite union of $\G_I$-orbits.
Any $H\in \Hcal_I$  defines a hyperplane $(H\cap I^\perp)/I$ in $I^\perp/I$ of hyperbolic signature, and so the corresponding hyperplane  in $N$ meets $C$. The hyperplane sections  of $C$ thus obtained, form a finite union of $\G(N)$ orbits. This implies that they define a decomposition $\Sigma_I$ of $C$ into (relatively open) rational locally-polyhedral cones in the sense that  any closed subcone of of the closure of $C$ that is spanned by finitely many rational vectors meets only finitely many members of $\Sigma$.  For any $\sigma\in \Sigma$ we denote by $N_\sigma$ its linear span in $N$ (a subspace defined over $\QQ$).  

Actually, $\Sigma$ is still incomplete as we have not included its `improper members' yet. These are the rays in the boundary of $C$ spanned by a rational isotropic vector. To such a ray $\sigma$  we are going to associate a vector space $N_\sigma$ which is in general \emph{not} its span.  In fact, our discussion of  $\Jcal^{(2)}_\Hcal$ already suggests which space to take here: first note that the line spanned  by  $\sigma$ corresponds to an isotropic plane  $I_\sigma$ in $V$ defined over $\QQ$ which contains $I$. Now recall that we denoted  by $J_{I_\sigma}$ the intersection of $I_\sigma^\perp$ and the $H\in\Hcal$ containing $I_\sigma$. The image of $J_{I_\sigma}$ in $N$ is then the intersection of the hyperplane $\sigma^\perp$ with all the images of the members of $\Hcal_I$ which contain $\sigma$. This will be our  $N_\sigma$.

The collection $\widehat{\Sigma}$ obtained by $\Sigma$ completing  as above (with the associated assignment  $\sigma\in\widehat{\Sigma}\mapsto N_\sigma$), are the data needed for a semi-toric embedding.  If $\Sigma$ were 
a decomposition into rational polyhedral cones (rather than only locally polyhedral cones), then we would have a torus embedding $T\subset T_\Sigma$ and if we then let then, following Mumford et al.\ \cite{amrt}, $\hat\Tcal^\Sigma$ be the interior of the closure of $\Tcal^\Sigma$ in $T_\Sigma$, then $\G(N)$ acts properly discretely  on $\hat\Tcal^\Sigma$ so that we can form the orbit space $\Tcal^\Sigma_{\G (N)}$ in the analytic category. This would give us the extension we want. But in this more general setting the situation is more delicate, for we will not always produce a locally compact space $\hat\Tcal^\Sigma$ (although it will be one after dividing out by the $\G(N)$-action). Here is how we proceed.
For any $\sigma\in \widehat{\Sigma}$ we denote by $\TT_\sigma\subset \TT$  the  subtorus $(N_\sigma\cap N_\ZZ)\otimes\CC^\times$ and by $T (\sigma)$ the quotient of the $\TT$-torsor $T$ by this subtorus and by $p_\sigma: T\to T_\sigma$ the projection.  We form the disjoint union
\[
\hat\Tcal^\Sigma:= \Tcal \sqcup\bigsqcup_{\sigma\in\widehat{\Sigma}} p_\sigma (\Tcal),
\]
a set which comes with an action of $\G (N)$. 
A basis for a $\G(N)$-invariant topology on $\hat\Tcal^\Sigma$ is defined as follows: 
given $\sigma\in\Sigma$, let  $\TT_\sigma^+\subset \TT_\sigma$ be the semigroup that is the preimage of $\sigma$ under the natural map $\Im : \TT_\sigma\to N_\sigma$. Then for any open subset $U\subset \Tcal$ that is convex,
and invariant under both $Z_{\G (N)}(\sigma)$ and the semi-group $\TT_\sigma^+$, we put 
\[
\hat U_\sigma:= U\sqcup\bigsqcup_{\tau\in\widehat\Sigma, \tau\le \sigma} p_\sigma (\U).
\]
and observe that such subsets define a topology on $\hat\Tcal^\Sigma$. The $\G (N)$-orbit space 
$\hat\Tcal^\Sigma_{\G (N)}$ is then locally compact Hausdorff and carries a structure as a normal analytic space for which the obvious maps $p_\sigma (\Tcal)\to  \hat\Tcal^\Sigma_{\G (N)}$ are analytic.

It is  possible to do this construction even before dividing out by the lattice $N_\ZZ$. We then can also carry along the $\CC^\times$-bundle $\LL^\times$ and  the construction as a set becomes even somewhat simpler to describe
(and more in spirit of the preceding):
let $\Jcal^{(I)}_\Hcal$ be the collection  of pairs  $(J, \sigma)$, where $J$ is an intersection of members of $\Hcal_I$ and either $\sigma\in \Sigma_I$ is such that its linear span is the image of $J$ in $I\otimes I^\perp/I$ or $\sigma$ is an improper ray such that for the corresponding isotropic plane $I_\sigma$ we have $J=J_{I_\sigma}$. Then
\[
\LL^\times \sqcup\bigsqcup_{(J, \sigma)\in\Jcal^{(I)}_\Hcal} \pi_J(\LL^\times).
\]
can be given a $\G_I$-invariant  topology, such its orbit space yields a $\CC^\times$-bundle over $\hat\Tcal^\Sigma_{\G (N)}$. \\

We have now defined $\Jcal_\Hcal=\cup_{i=0}^2\Jcal^{(i)}_\Hcal$, but  must still specify $\Jcal_\Hcal^\infty$: this will be the union of $\Jcal^{(1)}_\Hcal$ and the collection of $J_I$, where $I$ runs over the  isotropic $\QQ$-planes in $V$. Notice that  when $\Hcal$ is empty,  $\Jcal^\infty=\Jcal$ and equals  the collection of $I^\perp$,  where $I$ runs over the all the nonzero isotropic $\QQ$-subspaces in $V$. This then produces the  modification chain attached to our arithmetic arrangement.

\subsection*{A revisit of double covers of sextic curves}
As we noted the double covers of sextic curves do not produce all K3 surfaces of degree 2, even if we allow them to have rational double points. We are missing the K3 surfaces  that admit an elliptic fibration whose fibers have degree one withe respect to the polarization. They hide themselves behind the (closed semistable) orbit in $\PP(H^\sst)$ of triple conics. As Shah observed, these can be exposed by using twice the polarization: then our $\PP(U)$ is realized as the Veronese surface   $V\subset \PP^5$, where the sextic appears  as an intersection of $V$ with a cubic hypersurface in $\PP^5$. In particular, the triple conic appears as a triple hyperplane section of $V$.  If we do GIT with respect to the action of $\SL(6)$ on the Hilbert scheme of $\PP^5$ of such curves (of genus 10 and of degree 12), then it is not surprising that the latter is no longer semistable. The GIT compactification replace such curves by curves that are the intersection of a cubic hype surface with a projective cone $V_o$ whose base is a normal rational  curve in a $\PP^4$ (in other words $V$ degenerates into $V_o$). The double covers of $V_o$ which ramify over the union of  such a curve and the vertex then yield the missing K3 surfaces (with the projection from the vertex giving the elliptic fibration). We have thus accounted for $\Delta_\Hcal$.  The resulting compactification of the moduli space of K3 surfaces of degree 2 can be identified with our $\DD^\Hcal_\G$: it a small blowup of $\DD^\bb_\G$  with fibers of dimension at most $1$. As in the case of genus 3 curves, two distinct members  of $\Hcal$ do not meet in $\LL^\times$ and therefore no arrangement blowup is needed: the closure of $\Delta_\Hcal$ in $\DD^\Hcal_\G$ can be contracted to a variety of dimension $\le 2$ and this 
reproduces the GIT compactification of the moduli space of sextic plane curves.


\end{document}